\documentclass[a4paper,%
oneside,%
DIVcalc,%
]{scrbook}
\usepackage{diss}
\begin{document}
\frontmatter
%\address{Fachbereich Mathematik\\Universität Dortmund}
%\email{Boris.Hemkemeier@Math.Uni-Dortmund.DE}
\publishers{{\bf Versionshinweis:}\\ dies ist ein inhaltsgleicher, von arXiv
  formatierter Reprint einer Dissertation zur Erlangung des Grades
  eines\\ Doktors der 
  Naturwissenschaften\\[2cm]Im Original dem Fachbereich Mathematik der
  Universität 
  Dortmund vorgelegt im September~2003}
\title{Algorithmische Konstruktionen von Gittern}
\date{}
\author{Boris Hemkemeier}
\pagenumbering{roman}
 
\maketitle
\mainmatter
\pagenumbering{roman}
\thispagestyle{empty}
{}\hspace{4cm}

\vfill
 \begin{quote}\large
 I have  a great suspicion that for example Euler today would spend
 much more of his time on writing software because he spent so much of
 his time e.g., in efforts of calculating tables of moon
 positions. And I believe that Gau"s as well would spend much more time
 sitting in front of the screen. 

 \hfill\emph{Yuri I. Manin} \cite{manin98:_inter}
 \end{quote}
\vfill\vfill\vfill
\tableofcontents
\chapter[Einleitung]{Einleitung}
\label{sec:einleitung}
\pagenumbering{arabic}
Motivation dieser Arbeit ist ein Klassifikationsprojekt aus dem Gebiet
der ganzzahligen Gitter. Martin Kneser beschrieb 1957 in dem Artikel 
\emph{Klassenzahlen definiter quadratischer {F}ormen}
\cite{kneser57:_klass_formen} ein Exhaustionsverfahren f"ur die Klassifikation
vollst"andiger Geschlechter von Gittern. Diese Methode ist nicht unmittelbar
praktikabel f"ur ein automatisiertes Verfahren; sie l"a"st sich aber unter
algorithmischen Aspekten modifizieren, um auch Geschlechter mit gro"ser Klassenzahl
klassifizieren zu k"onnen.

Im Rahmen dieser Arbeit wurde daraus das Computerprogramm \texttt{tn}
entwickelt, um die Geschlechter vieler Gitter zu bestimmen, die bis auf wenige
Ausnahmen zuvor nicht vollst"andig bekannt waren.  Bei der Konstruktion neuer
Gitter mit Knesers Methode beh"alt man die Kontrolle "uber wichtige
Gitterinvarianten. Damit unterscheidet sich \texttt{tn} von lokalen
Optimierungs- oder Trainingsalgorithmen und hat eine gezielte Suche nach
Gittern mit herausragenden Eigenschaften wie Extremalit"at, Wur\-zel\-sys\-tem
oder hoher Packungsdichte erlaubt.

Die Algorithmik der Konstruktion ganzzahliger Gitter motivierte neue
Problemstellungen und erwies sich auch aus theoretischer Sicht als ein
ergiebiges Gebiet. Das scheinbar banale Problem der Zerlegung eines
Gitters in orthogonale Untergitter f"uhrte zur Entwicklung effizienter
Algorithmen zum einen f"ur das Zerlegungsproblem selbst, zum anderen für
die Erzeugung eines
Gitters aus einem gro"sen Erzeugendensystem. Minkowskis zweiter
Fundamentalsatz~\cite{minkowski96:_geomet_zahlen} aus der Geometrie
der Zahlen f"uhrt zu Absch"atzungen der Anzahl der ben"otigten
arithmetischen Operationen und zu einer allgemeinen, quantitativen
Aussage "uber die Gr"o"se minimaler Erzeugendensysteme von Gittern.

Der erste Teil der Arbeit leistet einen Beitrag zur konkreten Klassifikation
ganzzahliger Gitter. Die Fundierung und ein Teil der Ergebnisse wurden in der
gemeinsamen Arbeit \emph{Classification of integral lattices with large
  class number\/} \cite{scharlau98:_class} mit Rudolf Scharlau publiziert. In diesem Kapitel wird auf die origin"aren algorithmischen
Aspekte der $2$-Nachbarmethode fokussiert.
  Mit Hilfe von Knesers Nachbarverfahren gelang es
H.~Niemeier bereits 1973, das Geschlecht der geraden, unimodularen Gitter in Dimension $24$
zu bestimmen \cite{niemeier73:_defin_formen_dimen_diskr}.  Er konnte dabei die
starken Struktureigenschaften der Gitter in diesem Geschlecht ausnutzen, um
die Nachbargitter zu identifizieren und zu unterscheiden. In allgemeineren 
Situationen sind separierende Invarianten aber nur mit gro"sem algorithmischen
Aufwand zu berechnen. In diesem Sinne "au"sern sich Conway und Sloane
verhalten "uber den praktischen Nutzen von Knesers Methode in nichttrivialen
F"allen.
\begin{quotation}
  {"`However there is a geometric method used by Witt and Kneser which
    (after the work of Niemeier) is effective roughly until the sum of
    the dimension and the (determinant)$^{1/2}$ exceeds $24,$ beyond
    which point it seems that the forms are inherently
    unclassifiable."' \cite{bible}, S.~352~f.}
\end{quotation}
Diese Einsch"atzung erwies sich als zu pessimistisch, denn die Ergebnisse in
\cite{scharlau98:_class} und in dieser Arbeit zeigen, da"s
Klassifikationen deutlich jenseits dieser Schranke durch\-f"uhr\-bar sind.  

Zwei verschiedene Gitter $L$ und $L'$ hei"sen benachbart, wenn sie ein
gemeinsames Untergitter jeweils vom Index $2$ enthalten.  Bei geraden Gittern
mit ungerader Determinante l"a"st sich mit Knesers Nachbarmethode ein
Nachbargitter $L(v)$ aus $L$ und einem Nachbarvektor $v\in L\setminus 2L$ mit
Hilfe eines sogenannten Nachbarschrittes konstruieren. Da wir nur an den
Isomorphieklassen von Nachbargittern interessiert sind, l"a"st sich die
Konstruktion aller Nachbargitter bis auf Isomorphie auf eine endliche Auswahl
von Nachbarvektoren $v$ zur"uckf"uhren. Sie werden als geeignete
Repr"asentanten aus der endlichen Menge von Orbits der induzierten
Automorphismengruppe von $L$ auf $L/2L$ ausgew"ahlt. 

Das Computerprogramm \texttt{tn} basiert auf einer Implementation
dieser Methode. Daneben enth"alt es zahlreiche Analysefunktionen zur
Bestimmung von Invarianten der konstruierten Gitter.  Wegen der
Komplexit"at der Berechnungen wurde das Design und die Realisierung
von \texttt{tn} sehr stark an Performanceaspekten ausgerichtet und
deshalb nicht in einem existierenden Computeralgebrasystem
vorgenommen.\footnote{Seit 1999 existiert auch eine Implementierung
  der Basisfunktionalität der $2$-Nachbarmethode in dem
  Computeralgebrasystem MAGMA.}
  
Bei den enumerativen Klassifikationen konzentrierte sich das Interesse
auf die $\ell$-modularen Gitter mit $\ell=3,5,7,11$ in Dimensionen bis
zu $14$ und die Bestimmung ihrer modularen und extremalen Gitter. Ein
bemerkenswertes Ergebnis aus theoretischer Sicht ist ein Satz "uber
das Geschlecht $7^6.$ Die Klassifikation mit \texttt{tn} lieferte den
Beweis f"ur die Nichtexistenz eines vermuteten extremalen Gitters.

Das Programm \texttt{tn} kann zur Systematisierung der Ergebnisse eine
Vielzahl von Invarianten eines Gitters berechnen. Eine bisher wenig
beachtete Invariante ist die Fourier-Transformation der
L"angenfunktion eines Gitters auf $L/2L.$ Erste Vermutungen, da"s
diese Fourier-Transformierte nur einen kleinen Tr"ager besitzt,
bewahrheiteten sich nicht. Es scheinen jedoch im allgemeinen nur
relativ wenig verschiedene Fourierkoeffizienten in der transformierten
Funktion eine Rolle zu spielen. Mit \texttt{tn} lassen sich das
Spektrum und die Fourier-Transformierte der L"angenfunktion eines
Gitters für eine feste Basis bestimmen.

Die Nachbarschaftsmethode erweist sich auch in F"allen, in denen die Theorie
keinen Erfolg garantiert, als erstaunlich robust. Bei Gittern mit gerader
Determinante ist der oben beschriebene Nachbarschritt nicht immer m"oglich.
Dennoch läßt sich beispielsweise das (bereits bekannte) Geschlecht des
Barnes--Wall-Gitters \cite{MR95e:11073} in der Dimension $16$
berechnen.  Man kann beobachten, da"s in h"oheren Dimensionen der
Nachbarschaftsgraph, der gebildet wird aus den Isomorphieklassen von
Gittern als Knoten und den "Aquivalenzklassen von Nachbarvektoren als
Kanten, f"ur jedes Gitter so viele Kanten enth"alt, da"s das
Kollabieren einzelner Nachbarbildungen in der Regel nicht
ergebnisrelevant ist.

Das zweite Kapitel beginnt mit einer Diskussion des Zerlegbarkeitsproblems. Ein
Gitter hei"st zerlegbar, wenn es sich als direkte Summe nichttrivialer,
orthogonaler Untergitter schreiben l"a"st.
Die Zerlegung eines Gitters in unzerlegbare Gitter ist eindeutig
bis auf die Reihenfolge der Summanden. Der konstruktive Beweis dieses Satzes
stammt von Martin Kneser \cite{kneser54}.
Seine Konstruktion, eine Analogie zum Sieb von Erathostenes, dient prim"ar dem
Beweis der Behauptung, aber nicht der konkreten Durchf"uhrung. So
erweist sie sich unter dem Blickwinkel der Berechenbarkeit als nicht
ausreichend 
effizient. Wir zeigen in Satz~\ref{theorem:runtime_of_algorith_1}, wie ein
inkrementelles Konstruktionsverfahren statt der Benutzung des Siebverfahrens
zu praktikablen Laufzeiten f"uhrt, und geben eine Schranke f"ur die Anzahl der
algorithmischen Operationen des Gitterzerlegungsproblems an.

Das Fehlen eines Basiserg"anzungssatzes und eines Steinitz'schen
Basisaustauschsatzes in der Modulsituation verkompliziert die algorithmische
Behandlung der Basisberechnung eines Gitters au"serordentlich.  Das
Hauptergebnis in diesem Kapitel ist ein "Uberdeckungssatz f"ur
Erzeugendensysteme von Gittern, der aus dem zweiten Fundamentalsatz von
Minkowski  folgt.  Es wird gezeigt, da"s in
typischen Situationen der praktischen Berechnung von Gittern ein mininales
Erzeugendensystem nur aus wenig mehr Vektoren als dem Gitterrang besteht. Die
quantitative Aussage liefert Satz~\ref{satz:covering_theorem}: ist $S\subseteq
L$ ein minimales Erzeugendensystem eines Gitters $L$ auf dem $n$-dimensionalen
euklidischen Raum $E$  mit Vektoren, die in der
euklidischen Norm nicht l"anger als $B\in \BR$ sind, dann ist $|S|\leqslant
n
+ \log_2(n!(\frac{B}{\min L})^{n}).$ Mit Hilfe dieses Satzes wird gezeigt, da"s
das algorithmische Problem der Konstruktion von Gitterbasen aus gro"sen
Erzeugendensystemen durch das angegebene inkrementelle Verfahren effizient
l"osbar ist im Vergleich zu etablierten Methoden wie den Variationen des
LLL-Algorithmus und der Berechnung der Hermite-Normalform.  Tests mit
ganzzahligen Zufallsgittern unterstreichen die praktische Relevanz (siehe die
Laufzeitvergleiche auf Seite \pageref{fig:dim20_100}).

Das dritte Kapitel behandelt die Berechnung von Gitterquantizern. Dabei
handelt es sich um Approximationen, die aus der Informationstheorie motiviert
werden.  Ein Quantizer eines Gitters ordnet jedem Punkt im euklidischen Raum
den n"achstliegenden Gitterpunkt zu. Er ist also im wesentlichen durch die
Geometrie der einem Gitterpunkt n"achstliegenden Punkte, der sogenannten
Vorono"izelle bestimmt. Die Qualit"at eines Gitterquantizers aus
ingenieurwissenschaftlicher Sicht wird durch das dimensionslose zweite
normalisierte Tr"agheitsmoment der Vorono"izelle (Definition
\ref{def:traegheitsmoment}) bewertet. Dies ist ein skaliertes Integral "uber
die quadratische Vektorl"angenfunktion in der Vorono"izelle des Ursprungs. Die
Abstandsfunktion wird algorithmisch durch das Verfahren von Fincke und Pohst
auf Nebenklassen \cite{fincke-pohst} ausgewertet.\footnote{In der Praxis von
  Telekommunikationsanwendungen wie bei Analog-Digital-Wandlern spielt neben
  der G"ute des Tr"agheitsmomentes auch die Effizienz der Implementation der
  Approximation eine Rolle.} 
Wir beschreiben die bekannten, grundlegenden Eigenschaften von
Vorono"izellen und die approximative Berechnung ihrer Gitterquantizer mit
Monte-Carlo-Integration.  Ihre Berechnung ist ebenfalls in \texttt{tn}
implementiert. Alternativ lassen sich Gitterquantizer und Quantizer zu
Kugelpackungskonstruktionen aus nichtlinearen Codes (Typ A Konstruktion nach \cite{bible})
mit einem Programm f"ur das Computeralgebrasystem MAGMA 
berechnen. Gitterquantizer sind bis in die j"ungste Vergangenheit nicht
intensiv untersucht worden. In einer neueren Arbeit erzielten Agrell und
Eriksson \cite{agrell98:_optim_lattic_quant} eine Reihe von zum Teil
kontraintuitiven Ergebnissen über Gitterquantizer.

Für das Zustandekommen dieser Arbeit schulde ich vielen Personen
Dank. An erster Stelle ist hier Herr Prof.~Dr.~Rudolf Scharlau zu nennen. 
Seine Betreuung gab dieser Arbeit den Anfang und die Richtung. Er war
mir in all den Jahren ein aufmerksamer und geduldiger Gesprächspartner.
Die Diskussionen mit meinem Kollegen Frank Vallentin
haben mir viele Anknüpfungspunkte an die Informatik eröffnet. Bernd
Souvignier hat seine Implementationen der Berechnung von Erzeugern der
Automorphismengruppe eines Gitters und des  Isometrietests zweier
Gitter freundlicherweise diesem Projekt zur Verfügung
gestellt. 
Ohne die wunderbaren Entwicklungstools  emacs, gcc und gbd
aus dem GNU-Projekt von Richard Stallmans Free Software Foundation
wäre aus diesem Projekt  wahrscheinlich nicht ein so flexibles und
leistungsfähiges Programm wie  \texttt{tn} entstanden.
Die Hochschulrechenzentren der Universitäten 
Bielefeld, Dortmund und Duisburg haben mir Kapazitäten auf ihren
Supercomputern zur Verfügung gestellt.  
Der Besuch des  AT\&T Shannon Labs in New Jersey auf Einladung
von Neil Sloane 
 lenkte mein Interesse auf Gitterquantizer.
 Die  Universität Dortmund hat mir für die Anfertigung dieser
 Dissertation eine Assistentenstelle zur Verfügung gestellt.

Ganz besonders danke ich meiner Frau Dr.~Eva-Maria Kaffanke dafür, daß diese
Arbeit nach langer Zeit doch noch ihre Vollendung gefunden hat.

\chapter{Klassifikation ganzzahliger Gitter}\label{cha:tn}
\section{Einleitung}
\label{sec:tn_intro}

In diesem Kapitel wird das Nachbarverfahren zur Klassifikation
ganzzahliger Gitter mit gro"ser Klassenzahl beschrieben.
Die 
zugrundeliegende Methode wurde von M.~Kneser entwickelt
\cite{kneser57:_klass_formen}. H.-V.~Niemeier klassifizierte mit ihrer
Hilfe 1973 die geraden
unimodularen Gitter in Dimension $24$
\cite{niemeier73:_defin_formen_dimen_diskr}. 

Obwohl dieser Meilenstein in der Klassifikation der ganzzahligen
quadratischen Formen die Bedeutung von Knesers Methode beweist, wurde der
praktische Nutzen in "`generischeren"' Situationen als eher gering
eingesch"atzt (siehe \cite{bible}, S.~352).
In der gemeinsamen Arbeit mit
R.~Scharlau \emph{Classification of integral lattices with large class
  number} \cite{scharlau98:_class} wurde gezeigt, daß eine
algorithmisch optimierte Version des Nachbarverfahrens  eine Vielzahl
von Geschlechtern zug"anglich macht, die mit bisherigen Methoden nicht
klassifizierbar waren. Besonders für die Suche nach Gittern mit
bestimmten Eigenschaften kann ein Exhaustionsverfahren wie die
Nachbarmethode "außerst nützlich sein.
Zu dieser Arbeit ist das Computerprogramm
\texttt{tn} entstanden, dessen Funktionsweise im folgenden
erl"autert wird. In Abschnitt~\ref{sec:2neighbours} werden 
die theoretischen Grundlagen der Konstruktion von Gittern durch die
$2$-Nachbarmethode  zusammengefaßt.  In
Abschnitt~\ref{sec:program_tn} wird die Funktion von \texttt{tn}
genauer beschrieben und mit einigen markanten Resultaten in
Abschnitt~\ref{sec:tn_results} illustriert.\label{klass:unklare_numerierung} 

Knesers Verfahren garantiert unter gewissen Nebenbedingungen eine
vollst"andige Klassifikation von (Spinor-)Geschlechtern. Die $2$-Nachbarmethode
l"aßt sich aber auch mit Erfolg bei Gittern mit gerader Determinante
anwenden (Abschnitt \ref{sec:tn_even_det}), welche die Voraussetzungen
des Satzes von Kneser nicht erfüllen. Ein interessantes Beispiel ist die Klassifikation des Geschlechtes des
Barnes-Wall-Gitters (Abschnitt~\ref{sec:barnes_wall}). 

\section{Die Nachbarmethode}
\label{sec:2neighbours}
Im folgenden seien alle Gitter ganzzahlig mit ungerader Determinante
auf dem euklidischen $n$-dimensionalen Raum $(V,(-,-)),$ sofern nicht
etwas anderes vermerkt ist. 

\begin{definition}\label{def:2nachbar}
  Zwei Gitter $L$ und $L'$ werden Nachbarn genannt, falls
  \begin{displaymath}
    L/(L\cap L') \cong \BZ/2\BZ \cong L'/(L\cap L').
  \end{displaymath}
\end{definition}
Aus dieser Definition folgt unmittelbar das 
\begin{lem}
  Zwei benachbarte Gitter $L$ und $L'$ besitzen dieselbe Determinante.
\end{lem}
\begin{proof}
  Nach der Determinanten-Index-Formel gilt $$\det L = [L:L\cap L']^2 \cdot
  \det L\cap L' = [L':L\cap L']^2 \cdot \det L\cap L' = \det
  L'.\qquad\popQED$$ 
\renewcommand{\qedsymbol}{}
\end{proof}

Die effektive Konstruierbarkeit von Nachbarn kl"art das 
\begin{lem}\label{lem:construction_of_neighbours}\rule{0cm}{1mm}
  \begin{enumerate}
  \item 
    F"ur ein gerades Gitter $L$ und einen Vektor $v\in L\setminus 2L$ ist das
    Gitter
    \begin{displaymath}
      L(v) := L_v+ \frac{1}{2} \BZ v \quad \text{mit}\quad L_v:= \{x\in L\mid
      (x,v)\in 2\BZ \}
    \end{displaymath}
    ein Nachbar von $L.$ Man nennt $L(v)$ den Nachbarn von $L$
    bez"uglich des Nachbarvektors $v.$
  \item Wenn die Gitter $L$ und $L'$ Nachbarn sind, dann gilt
    $L(2w)=L'$ f"ur jedes $w\in L'\setminus L.$
  \end{enumerate}
\end{lem}
\begin{proof}
  Zu i) Wir zeigen zun"achst $|L/L_v|=2.$ Die auf $\overline{L}:=L/2L$
    induzierte Bilinearform $(\overline{v},\overline{w}):=\overline{(v,w)}$
    besitzt Diskriminante $1$ und ist insbesondere regul"ar, weil die
    Diskriminante von $(-,-)$ ungerade ist. Wir nehmen an, da"s $L=L_v.$ Dann
    gilt $(v,w)\in 2\BZ$ für alle $w\in L.$ Somit ist
    $(\overline{v},\overline{w})=0$ f"ur alle $\overline{w}\in \overline{L},$
    also ist $\overline{v}=0$ und damit $v\in 2L,$ was im Widerspruch zur
    Voraussetzung steht.

    Aus der Exaktheit der Sequenz 
    \begin{displaymath}
      0 \longrightarrow L_v \longrightarrow L \xrightarrow{x\mapsto
      (\overline{v}, \overline{x})} \BZ/2\BZ
      \longrightarrow 0
    \end{displaymath}
    folgt, da"s $|L/L_v|=2.$
    
    Dann gilt offenbar auch $L/(L\cap L(v))=2.$ Da $L$ gerade ist, folgt, daß
    $v\in L_v,$ aber $\frac{1}{2}v\notin L_v,$ also ist $L(v)/(L\cap
    L(v))\cong \BZ/2\BZ$
    und damit sind $L$ und $L(v)$ Nachbarn.
  
    Zu ii). Es sei $w\in L'\setminus L$ für die Nachbarn $L$ und $L'.$
    Dann folgt 
    aus $|L'/(L\cap L')|=2,$ daß $2w\in L\setminus{2}L.$
    Wie oben ergibt sich aus
    \begin{displaymath}
      0 \longrightarrow L_{2w} \longrightarrow L \xrightarrow{x\mapsto
        (\overline{2w}, \overline{x})} \BZ/2\BZ
      \longrightarrow 0,
    \end{displaymath}
    da"s $|L/L_{2w}|=2.$ Weil $L'$ ganzzahlig ist, gilt f"ur jedes $x\in L',$
    da"s $(x,2w)=2(x,w)\in 2\BZ,$ also insbesondere $L\cap L'\subseteq
    L_{2w}.$ Da $L/L\cap L' \cong \BZ/2\BZ \cong L/L_{2w},$ mu"s $|L_{2w}/L\cap
    L'|=1$ sein, also ist $L_{2w}=L\cap L'.$ Aus $w\in L'\setminus L$ und
    $|L'/L_{2w}|=2$ erhalten wir schlie"slich $L'=L_{2w}+\frac{1}{2}\BZ (2w) =
    L(2w).$
\end{proof}
\begin{remark}
  $L(v)$ ist nicht notwendigerweise ganzzahlig. Dieser  Sachverhalt
  wird   im folgenden Lemma genauer untersucht.
\end{remark}

\section{Die Konstruktion des Nachbarschaftsgraphen}
\label{sec:konstruiere_nachbarschaftsgraph}

\begin{lem}\label{lem:integral_even}
Sei $L$ ein gerades Gitter und $v\in L\setminus 2L.$ Dann ist
\begin{enumerate}
\item $L(v)$ ein ganzzahliges Gitter genau dann, wenn $(v,v)\in 4\BZ,$
\item $L(v)$ ist ein gerades Gitter genau dann, wenn $(v,v)\in 8\BZ.$ 
\end{enumerate}
\end{lem}
\begin{proof}{}\hspace{4cm}{}

  \begin{enumerate}
  \item Sei $(v,v)\in 4\BZ.$ Mit $L$ ist auch $L_v$ ganzzahlig. Aus 
    $(L_v,\frac{1}{2}v)= \frac{1}{2}(L_v,v) \subseteq \BZ$ und
    $(\frac{1}{2}v,\frac{1}{2}v)= \frac{1}{4}(v,v)\subseteq \frac{1}{4}
    4\,\BZ=\BZ,$ folgt, daß  auch $L(v)$ ganzzahlig ist. Die Umkehrung folgt
    unmittelbar aus der Ganzzahligkeit von $L$ und $(v,v)\in \BZ
    \Leftrightarrow (\frac{1}{2}v,\frac{1}{2}v)\in 4\BZ.$
  \item  Sei $(v,v)\in 8\BZ.$ Mit $L$ ist auch $L_v$ gerade. Dann ist wegen
    $(\frac{1}{2}v,\frac{1}{2}v) \in \frac{1}{4}\,8\BZ=2\BZ$ auch $L(v)$
    gerade. Die Umkehrung folgt
    unmittelbar aus der Ganzzahligkeit von $L$ und $(v,v)\in 2\BZ
    \Leftrightarrow (\frac{1}{2}v,\frac{1}{2}v)\in 8\BZ.$
  \end{enumerate}
\end{proof}

Die beiden folgenden Lemmata erlauben es, die Konstruktion aller
Nachbarn eines Gitters bis auf Isomorphie
auf ein endliches Problem zur"uckzuf"uhren.
\begin{lem}
  Es seien $w,w'\in L\setminus 2L.$ Dann gilt $L(w)=L(w'),$ falls $w-w'\in
  2L_w.$
\end{lem}
\begin{proof}
  Es sei $z:=\frac{w-w'}{2}\in L_w.$ Dann folgt aus  $\overline{(x,w)}=
  \overline{(x,w)}-\overline{(x,2z)}= \overline{(x,w-2z)}= \overline{(x,w')}$
  f"ur alle $x\in L,$ daß $L_w=L_{w'}.$ Da $2L\subseteq L_w=L_{w'},$ folgt
  \begin{displaymath}
    L(w)=L_w+\frac{1}{2}\BZ w= L_w+\frac{1}{2}\BZ (w'+2z)=L_{w'}+\frac{1}{2}\BZ
    w'= L(w'). 
  \end{displaymath}
\end{proof}
\begin{lem}\label{lem:neighbour_even_odd}
  Sei $L$ ein ganzzahliges, gerades Gitter und $v\in L\setminus 2L.$
  Dann existiert ein $y\in
  L$ mit $(v,y)\notin 2\BZ$ und f"ur alle $w\in \overline{v}\in L/2L$ gilt
  entweder $L(w)=L(v)$ oder $L(w)=L(v+2y).$ Weiter ist eines dieser Gitter
  gerade, das andere ungerade.
\end{lem}
\begin{proof}
  Da $\overline{v}\neq 0$ erfüllt sogar für jede gegebene Basis von
  $L$ einer der Basisvektoren die von $y$ geforderte Eigenschaft. Ansonsten w"are die 
  induzierte Form $\overline{(-,-)}$ identisch $0$ und bes"aße damit
  nicht die Determinante $1.$
 
 Falls $(v,v)\equiv 0
  \pmod{8},$ so ist $L(v)$ ein gerades Gitter und $L(v+2y)$ wegen
  \begin{displaymath}
    (v+2y,v+2y)=(v,v)+4(v,y)+4(y,y)\equiv 4\pmod{8}
  \end{displaymath}
  ein ungerades Gitter. Insbesondere sind $L(v)$ und $L(v+2y)$ verschieden.
  Also gilt nach dem vorangehenden Lemma f"ur alle $w\in \overline{v}$ wegen
  $[2L:2L_v]=2$ entweder $L(w)=L(v)$ oder $L(w)=L(v+2y).$
\end{proof}

Damit braucht man lediglich $2^{n}-1$ Repr"asentanten der Klassen von
$L/2L\cong \BF_2^n$ als Nachbarvektoren zu ber"ucksichtigen, um alle
Nachbargitter von $L$ zu konstruieren. Dies reduziert die Komplexit"at der
Auswahl der Nachbarvektoren auf ein endliches Problem.

Dieses Argument beweist das
\begin{kor}
  Ein Gitter besitzt bis auf Isometrie h"ochstens  $2^n-1$ Nachbargitter.
\end{kor}

Das n"achste Lemma zeigt, da"s Nachbarvektoren, die in demselben Orbit unter
der Automorphismengruppe des Gitters liegen, isometrische Nachbarn erzeugen.
\begin{lem}
  Ist $\sigma\in O(L)$ und $w\in L\setminus 2L,$ dann ist
  $L(w\sigma)=L(w)\,\sigma.$ 
  \begin{proof}
    F"ur $x\in L$  gilt
    \begin{displaymath}
      x\in L_{w\sigma} \Leftrightarrow
      (x,w\sigma)\in 2\BZ \Leftrightarrow
      (x\,\sigma^{-1},w)\in 2\BZ\Leftrightarrow
      x\,\sigma^{-1}\in L_w \Leftrightarrow
      x\in L_w\, \sigma.
    \end{displaymath}
  \end{proof}
\end{lem}
\begin{remark}\label{rem:nachbarn}
  Aufgrund des vorhergehenden Lemmas können wir uns bei der Auswahl von
  Nachbarvektoren auf diejenigen Gittervektoren beschr"anken, die
  einerseits in   verschiedenen Orbits unter der
  Automorphismengruppe $O(L)$ liegen und die andererseits in
  verschiedenen Klassen modulo $2L$ liegen. Das 
  kommutative Diagramm in Abbildung~\ref{fig:orbits_unter_O_von_L} 
  veranschaulicht diesen Sachverhalt. Alle enthaltenen Abbildungen
  sind surjektiv. Wegen der
  Kommutativit"at des Diagramms 
  kann man die Berechnung
  aller Nachbargitter (bis auf Isometrie) von $L$ vereinfachen, indem
  man zun"achst zu einem Repr"asentantensystem von $L/2L\cong \BF_2^n$
  "ubergeht und darauf die induzierte Automorphismengruppe
  $\overline{O(L)}$ in kanonischer Weise
  operieren l"aßt. Hierbei bildet man Vektoren aus $L$ nach $F_2^n$ ab, verm"oge
  $\sum_{i=1}^n \alpha_i b_i \mapsto
  (\overline{\alpha_1},\cdots,\overline{\alpha_n})$ f"ur eine
  beliebige, aber feste Basis
  $\{b_1,\cdots,b_n\}$ von~$L.$
  \begin{figure}[htbp]
    \begin{center}
      
      \begin{displaymath}
        \begin{CD}
          L @>O(L)>> \{v\,O(L)\mid v\in L\}\\
          @VVV @VVV\\
          {\BF_2}^n@>\overline{O(L)}>>\{\overline{v}\,\overline{O(L)}\mid
          \overline{v}\in \BF_2^n\}
        \end{CD}
      \end{displaymath}
      \caption{Die Orbits unter $\overline{O(L)}$}
      \label{fig:orbits_unter_O_von_L}
    \end{center}
  \end{figure}

Die Konstruktion aller Isomorphieklassen von Nachbargittern des
Gitters $L$ ist folgendermaßen möglich:
\begin{enumerate}
\item W"ahle eine beliebige, von nun an fixierte Basis von $L.$
\item Identifiziere $L/2L$ mit $\BF_2^n$ 
  und berechne alle Orbits von $\overline{O(L)}$ auf $\BF_2^n$ bezüglich
  der gew"ahlten Basis. 
\item W"ahle aus jedem Orbit in $\BF_2^n$ unter der
  Operation von $\overline{O(L)}$ einen Vektor $\overline{v}\in
  \BF_2^n$ und suche ein beliebiges Urbild $v\in L$ auf.
\item Bilde den geraden Nachbarn von $L$ bezüglich $v,$ sofern er
  überhaupt existiert.
\end{enumerate}
Die Isomorphieklasse des Nachbargitters h"angt nicht vom gew"ahlten
Repr"asentanten ab, deshalb sprechen wir auch gelegentlich von einer
Nachbarbildung "`bezüglich der Klasse $\overline{v}$"' anstatt
"`bezüglich des Vektors $v$\grqq.
\end{remark} 
Mit der iterierten Bildung von Nachbarn eines gew"ahlten Gitters werden  
 weitere Gitter konstruiert. Bei geraden Gittern besteht
diese Menge aus einer Vereinigung von echten Spinorgeschlechtern, bei allen
 Beispielen und Resultaten in diesem Kapitel sogar aus einem ganzen
Geschlecht (siehe~\cite{scharlau98:_class}, S.~742; \cite{o'meara71:_introd},
\S 102 und \cite{kneser56:_klass_formen}, Satz~2 bis
Satz~5). Es bietet sich an, die Nachbarschaft von Gittern
als Inzidenzrelation in einem Graphen zu beschreiben. Dieser Graph ist
ungerichtet, weil nach
Lemma~\ref{lem:construction_of_neighbours}~ii)~{}Nachbarschaft eine
symmetrische Relation ist.

Da wir im folgenden an einer Klassifikation dieser Geschlechter
interessiert sind, definieren wir den Nachbarschaftsgraphen nicht 
auf der Menge aller Gitter, sondern nur innerhalb eines Geschlechts
und nur bis auf Isomorphie der Gitter.

\begin{definition}
  Für ein Gitter $L$ sei $[L]$ die Klasse aller zu $L$ isometrischen Gitter.
  Dann nennen wir den ungerichteten Graphen $(V,E)$ mit $V=\{[L]\mid
  L\in \mathcal{G} \}$ und 
  \begin{displaymath}
    E=\{\ \{[L_1],[L_2]\}\mid L_1,L_2\in \mathcal{G},
    \quad L_1 \text{ und } L_2 
    \text{ sind Nachbargitter}\} 
  \end{displaymath}
den \emph{Nachbarschaftsgraphen\/} von  $\mathcal{G}.$ 
\end{definition}

\section{Konstruktion einer Basis eines Nachbarn}
\label{sec:basiskonstruktion}

Sei $\{b_1,\cdots,b_n\}$ eine Basis  des geraden Gitters $L$ und $v\in
L\setminus 2L.$ In Algorithmus~\ref{algo:basis_of_even_neighbour} wird eine Basis des geraden Nachbarn
bez"uglich $\overline{v}$ konstruiert. Im Anschlu"s zeigen wir die Korrektheit
des Verfahrens.

\begin{algorithm}[h]
\caption{Eine Basis des geraden Nachbarn von  $L$ bzgl.~$v$}
\label{algo:basis_of_even_neighbour}
\begin{algorithmic}
  \STATE \emph{Input:} Eine Basis $\{b_1,\cdots,b_n\}$ von $L$ und
  $v=(v_1,\cdots, v_n)\in
  L\setminus 2L$ mit $(v,v)\in 4\BZ.$

  \STATE {\emph{Output:} Eine Basis $\{b_1',\cdots,b_n'\}$ des geraden Nachbarn
  von $L$ bez"uglich $v,$ falls dieser existiert, ansonsten eine Fehlermeldung} 
  \IF{$(v,v)\equiv 2\pmod{8}$ oder $(v,v)\equiv 6\pmod{8}$}
  \STATE{return FAILURE}
  \ENDIF

  \STATE {\sl // 1.\ Initialisierung.}
  \IF{$(v,v)\equiv 4\pmod{8}$}
  \STATE{w"ahle $i$ mit $(v+b_i,v+b_i)\equiv 0\pmod{8}$}
  \STATE{$v\leftarrow v+b_i$}
  \ENDIF
  \STATE{\sl // 2. Normalisiere $v.$}
  \STATE {$k\leftarrow \min \{i\,|\,v_i\notin 2\BZ\}$}
  \STATE{w"ahle $x$ mit $x\cdot v_k\equiv 1 \pmod 8 $}
  \STATE{$v\leftarrow x\cdot v$}
  \STATE {\sl // 3.\ Konstruktion der Basisvektoren.}
  \STATE{w"ahle $m\neq k$ mit $(v,b_m)\notin 2\BZ$}
  \STATE {$b_k'\leftarrow \frac{1}{2} v, \quad b_m'\leftarrow 2b_m$}
  \FOR{$i=1,\cdots,n,\quad k\neq i\neq m$}
  \IF{$(b_i,v)\in 2\BZ$}
  \STATE{$b_i'\leftarrow b_i$}
  \ELSE
  \STATE{$b_i'\leftarrow b_i+b_m$}
  \ENDIF \COMMENT   {  // also $(b_i',v)\in 2\BZ$ }
  \ENDFOR
  \STATE{return $b'$}
\end{algorithmic}
\end{algorithm}

In der Initialisierung im ersten Schritt wird $v,$ wenn m"oglich, auf
eine durch $8$ teilbare Norm gebracht.
Ist $(v,v)\equiv 0\pmod{8},$ so bilden wir den Nachbarn $L(v).$ Ist
$(v,v)\equiv 4\pmod{8},$ so w"ahlen wir nach 
Lemma~\ref{lem:neighbour_even_odd} ein $b_i$ mit
$(v,b_i)\notin 2\BZ$ und bilden $L(v+2b_i).$ 

Wir suchen den ersten ungeraden Koeffizienten $v_k$ in $v$ auf und
multiplizieren $v$ mit einem geeigneten Skalar $x$, so da"s $x\cdot v_k\equiv
1\pmod 8.$ In der Implementation wird $x$ durch einen einfachen Table--Lookup
bestimmt.
\begin{displaymath}
  x = \left\{
    \begin{array}{ll}
      3,& \text{ falls } v_k=3\\
      5,& \text{ falls } v_k=5\\
      7,& \text{ falls } v_k=7\\
    \end{array}
    \right..
\end{displaymath}
Da alle skalaren Faktoren ungerade sind, gilt immer noch $L(v)=L(xv).$
Anschlie"send wird ein $m\neq k$ bestimmt mit $(v,b_m)\notin 2\BZ.$
Angenommen, ein solches $m$ existiere nicht. Sei
$M=\Zspan{\{b_i\}_{i\neq k}}.$ Dann liegt nach der Annahme $\overline{v}\in \overline{M}^\perp,$
also ergibt sich beim "Ubergang zum orthogonalen Komplement
\begin{math}
  \overline{M}\subseteq \overline{v}^\perp.
\end{math}
Aus Dimensionsgr"unden gilt
\begin{math}
  \overline{M} = \overline{v}^\perp.
\end{math}
Da $(v,v)\in 2\BZ,$ folgt
\begin{math}
 v\in  \overline{v}^\perp= \overline{M}.
\end{math}
Wegen $\overline{v_k}=\overline{1}$ folgt $v\notin M,$ was ein Widerspruch ist.

Aus der 
Konstruktion folgt, daß $b_i' \in L_v,$ für $i\neq k$  und daß $b_k'\in L(v),$ also gilt 
\begin{math}
\Zspan{b_1',\cdots,b_n'}\subseteq L(v).
\end{math}
Die "Ubergangsmatrix von 
$\{b_1,\cdots,b_n\}$ zu $\{b_1',\cdots,b_n'\}$  besitzt folgende Gestalt
 (hierbei nehmen wir ohne Beschr"ankung der Allgemeinheit an, da"s $k=1$ und
$m=n$ ist): 

\begin{displaymath}
  \begin{pmatrix} \frac{1}{2} & 0 & 0 &\cdot  &\cdot &\cdot & 0 \cr
          \ast  & 1 & 0 &\cdot  &\cdot &\cdot & 0 \cr
          \ast  & 0 & 1 &0  &\cdot &\cdot &0  \cr
            \cdot   & \cdot & 0 & 1 &\cdot &\cdot &\cdot  \cr
            \cdot   & \cdot & \cdot & \cdot &\cdot &\cdot &\cdot  \cr
          \ast  & 0  & 0   & \cdot &0 &1 & 0 \cr
          \ast  &\ast& \ast& \cdot &\ast &\ast &2 \cr                
\end{pmatrix}
\end{displaymath}
Die Determinante dieser "Ubergangsmatrix ist gleich $1;$ damit folgt
aus dem Vergleich der Determinanten
\begin{math}
  \det \Zspan{b_1',\cdots,b_n'} = \det L = \det L(v),
\end{math}
da"s
$\Zspan{b_1',\cdots,b_n'}=L(v).$ 

\subsection{Ungerade Gitter und der duale Nachbarschaftsgraph}
\label{sec:dual_graph}

Nach Lemma \ref{lem:neighbour_even_odd} gehört zu jedem geraden Nachbarn $L'$
eines geraden Gitters $L$ genau ein ungerader Nachbar $L''.$ Dieser
Sachverhalt l"aßt sich auch von einem elementaren algebraischen Standpunkt aus
betrachten. Sei $v\in L\setminus 2L$ mit $L(v)=L'$ und $y\in L$ mit
$(v,y)\notin 2\BZ$ und $(v,v)\equiv 0\pmod{4},$ d.h., $L(v+2y)=L''.$
Wir setzen $M=L+\frac{1}{2}\BZ v$ und 
erhalten mit $M/L_v$ die elementarabelsche Gruppe der Ordnung $4$ mit den drei
nichttrivialen Elementen $L/L_v, L'/L_v, L''/L_v.$ Mit $L$ und $v$ ist die
Gruppe $M/L_v$ und damit das Paar aus dem ungeraden und aus dem
geraden Nachbarn von $L$ 
 eindeutig bestimmt. Dies
liefert uns eine (wohlbekannte) Interpretation des Nachbarschaftsgraphen. Eine
Kante im Graphen kennzeichnet eine Nachbarschaft von geraden Gittern. Dieser
Kante ordnen wir das entsprechende ungerade Gitter zu. Der duale Graph ist
also eine Teilmenge des Nachbarschaftsgraphen für ungerade Gitter. Bei
ungeraden Gittern gibt es im allgemeinen keine Korrespondenz von geraden und
ungeraden Nachbarn, wie man leicht an dem Fall sieht, daß $v$ in dem
orthogonalen Komplement  eines ungeraden Teilgitters $N\subset L$ liegt. Dann
ist auch $N\subseteq L_v,$ also sind beide Nachbarn von $L$ bezüglich $v$
ungerade. 

In Spezialf"allen kann der Übergang zum dualen Graphen "außerst
nützlich sein. R.E.~Borcherds benutzt in \cite{borcherds:_dimen_odd_unimod_lattic}
den Nachbarschaftsgraphen der Niemeier-Gitter, um die ungeraden unimodularen
Gitter in Dimension $24$ zu klassifizieren.

\section{Das Programm \texttt{tn}}
\label{sec:program_tn}
Die konkrete Berechnung des Nachbarschaftsgraphen erweist sich als eine
komplexe Aufgabe. Dies liegt an der Vielzahl der involvierten
Teilprobleme und verwendeten Zahlsysteme. Im Programm selbst werden
drei verschiedene Arithmetiken eingesetzt. Der aufwendigste
Teil in den Berechnungen von \texttt{tn} besteht in der Auswahl der geeigneten
Vektoren zur Nachbarbildung. Das Implementationsziel von \texttt{tn} ist, alle
Nachbarbildungen berechnen zu können, solange es "uberhaupt m"oglich ist, die
Automorphismengruppen der beteiligten Gitter zu berechnen und die Gitter auf
Isometrie zu pr"ufen.

Das Programm \texttt{tn} behandelt nur ganzzahlige Gitter, so daß die
meisten Berechnungen in hardwarenaher $32$-Bit Integerarithmetik erfolgen
können.

Einige Routinen verlangen die Berechnung der "`kurzen Vektoren"' eines
Gitters; damit ist eine Auf\/listung oder Abz"ahlung aller
Gittervektoren bis zu einer bestimmten L"ange, meist der L"ange des
größten Basisvektors, gemeint. Dazu benutzen wir den Algorithmus von
U.~Fincke und M.~Pohst \cite{fincke-pohst} mit den Erg"anzungen von
H.~Cohen (Algorithm~2.7.7 in
\cite{cohen96:_cours_comput_algeb_number_theory}, S.~105) in einer
Implementation mit Doublearithmetik.

Die Berechnungen zur Auswahl der Vektoren zur Nachbarbildung werden,
wie in Bemerkung~\ref{rem:nachbarn} beschrieben, haupts"achlich  im Vektorraum
$L/2L\cong \BF_2^n$ vorgenommen. Für diesen Zahlbereich wurde  eine eigene
Implementation in Bin"ararithmetik implementiert, die in Abschnitt~\ref{sec:binary_arithmetic} vorgestellt wird.
\subsection{Die Berechnung der Automorphismengruppe eines Gitters}
\label{sec:isoautom}

Für die Berechnung von  Automorphismengruppen und Isometrieen von Gitteren
werden in \texttt{tn} die Programme \texttt{autom} und \texttt{isom}
verwendet, die 
uns Bernd Souvignier freundlicherweise zur Verfügung stellte.\footnote{Bernd
  Souvignier implementierte sp"ater neben \texttt{isom} und
  \texttt{autom} die $2$-Nachbarmethode in dem Computeralgebrasystem MAGMA.} Sie sind
Implementation der von B.~Souvignier und W.~Plesken in
\cite{plesken97:_comput_isomet_lattic} beschriebenen Verbesserungen der
generischen Sims-Stabilisatorkettenmethode in Gittern. 

Die Routinen \texttt{autom} und \texttt{isom} sind vom Rechenzeitbedarf her
die aufwendigsten Programmteile von \texttt{tn}. 

\subsection{Rechnen in $L/2L$}\label{sec:binary_arithmetic}
Da $\overline{L}= L/2L\cong \BF_2^n,$ liegt es nahe, die Rechnungen in
$\overline{L}$ in hardwarenaher Bin"ararithmetik zu
implementieren. Die Berechnung aller Orbits unter $\overline{O(L)}$
erfolgt durch Anwendung aller 
von \texttt{autom} berechneten und anschließend induzierten
Erzeuger. Aus diesem Grund ist die Effizienz dieser Routinen im
Hinblick auf Geschwindigkeit und Speicherplatzbedarf in der Praxis ein
kritischer Erfolgsfaktor des Programms \texttt{tn}. Die Implementation erfolgte
in der Programmiersprache "`C"' \cite{kernighan_ritchie}. 

Vektoren aus $\BF_2^n$ werden mit dem Datentyp \texttt{unsigned int}
identifiziert, der in den meisten Prozessorarchitekturen in 
einem 32-Bit Wort implementiert ist. In Dimensionen jenseits von $32$ sind aus
mehreren Gründen keine kompletten Klas\-si\-fi\-ka\-ti\-on\-en zu erwarten. Auf der
einen Seite dauern Automorphismengruppenberechnungen und Isometrietests
bereits sehr lange, auf der anderen Seite ist die Klassenzahl bei
Geschlechtern in diesen Dimensionen sehr gro"s. Beispielsweise wei"s man durch
Anwendung der Ma"sformel (s.u.), da"s es mehr als achtzig Millionen
unimodulare gerade Gitter in Dimension $32$ gibt.  

Bei der Berechnung der Orbits von $\overline{O(L)}$ auf $\overline{L}$
verwenden wir zwei verschiedene Kodierungen f"ur Vektoren in $\BF^n,$ je nach
Optimierungsziel, Rechenzeitbedarf und Speicherplatzbedarf. Bei der
Matrix-Vektor-Multiplikation wird ein Vektor von $\overline{L}$ als Bitfolge in einem Objekt
vom Typ \texttt{unsigned int} kodiert. Eine Matrix besteht aus einem Array von
$n$ solcher Objekte. F"ur die Markierung der Orbits von $\overline{O(L)}$
verwenden wir im Vektorraum $\overline{L}$ ein Kompressionsverfahren, bei dem
jeder Vektor in $\overline{L}$ nur durch ein einziges Bit
repr"asentieren wird, und die Kodierung und Dekodierung sehr effizient erfolgt.
 Die
folgenden beiden Beispiele illustrieren die benutzten Techniken bei der
Matrix-Vektor Multiplikation in $\overline{L}$
(\ref{subsub:mv-multiplikation}) und der Bestimmung eines kanonischen
Repr"asentanten eines Orbits von $\overline{O(L)}.$ Der Quelltext von
\texttt{tn} ist hier
f"ur den Zweck der  Darstellung aufbereitet worden. Fehlerbehandlung und der
Umgang mit Spezialf"allen wie $\dim L \leqslant 3$ f"uhren zu  etwas
komplizierterem Code.

\subsubsection{Matrix-Vektor-Multiplikation in
  $\overline{L}$}\label{subsub:mv-multiplikation} 

Als Beispiel f"ur die Effizienz der Arithmetik in $\overline{L}$ geben wir
hier die verwendete Routine für eine Matrix-Vektor-Multiplikation an.

\lstset{escapechar=\%}
\begin{lstlisting}[caption=%
{Effiziente Matrix-Vektor Multiplikation in $\overline{L}$}]{}
/* Matrix-Vektor Multiplikation in L/2L */

typedef unsigned int vector_2;

/* A definition for matrices over  GF(2). */
typedef
struct struct_matrix_2
{
  int dim;       /* dimension */
  vector_2 *v;   /* rows of the matrix */
}
*matrix_2;

void
matrix_vector_2 (vector_2 *imv, const matrix_2 m,
                 const vector_2 v)
     /* matrix_vector_2 computes the image of */
     /* the column vector v under the */
     /* operation of the matrix m and puts the */
     /* result in imv. */
{
  int i, j;
  vector_2 tmp;

  for (i = 0, *imv = 0; i < m->dim; i++){
     /* Compute (m->v[i],v) and set result
        in i-th bit of im */
     tmp = m->v[i] & v; 
     for (j = __WORDSIZE / 2; j >= 1; j /= 2) 
        tmp ^= tmp >> j; %\label{lst:fold_tmp}%
     *imv |= (tmp & 1) << ((m->dim - 1) - i); %\label{lst:ret_im}%
    }
}
\end{lstlisting}

Nachdem der Bildvektor \texttt{imv} mit dem Nullvektor initialisiert ist, wird
in der \texttt{for}-Schleife jeweils die $i$-te Komponente von \texttt{imv}
ermittelt. Hierzu berechnet man das Skalarprodukt der $i$-ten Zeile
\texttt{m->v[i]} der Matrix \texttt{m} mit dem Vektor \texttt{v} und weist 
es dem $i$-ten Bit in \texttt{im} zu. Man setzt $j$ auf die halbe Anzahl der
Bits in einem \texttt{vector\_2}. 

Die Zeile~\ref{lst:fold_tmp} wird in der Schleife durchlaufen, um durch
mehrfache Faltung mit einer XOR-Operation das Ergebnis im $0$-ten Bit von
\texttt{tmp} zu erhalten. Die oberen Bits spielen f"ur das Ergebnis dann keine
Rolle mehr. Das $0$-te Bit wird in Zeile~\ref{lst:ret_im} an die $(n-i)$-te
Stelle in \texttt{im} kopiert. Um bereits gesetzte Bits nicht zu l"oschen,
erfolgt dies mit einer logischen OR-Operation auf \texttt{im.}

\subsubsection{Orbits unter Gruppenoperation}

Die allgemeine Standardmethode, einen Orbit unter einer
Gruppenoperation zu bestimmen, ist der elegante \texttt{union-find}
Algorithmus mit Pfadkompression \cite{tarjan84:_union_find}.  In
unserer Situation können wir zus"atzlich ausnutzen, da"s die Orbits
aller Vektoren in $\overline{L}$ bestimmt werden müssen. Dies erlaubt
ein effizienteres Vorgehen sowohl in Bezug auf Speicherplatz als auch
auf Rechenzeit.  Dabei identifizieren wir $\overline{L}$ mit den
Koordinatenvektoren, die durch eine feste Basis von $L$ gegeben sind
und kodieren jeden Koordinatenvektor auf ein bestimmtes Bit.  Diese
Bits bilden eine Folge in einem zusammenh"angenden Speicherbereich.
Dadurch wird auf $\overline{L}$ eine totale Ordnung induziert, die bei
unserer Kodierung mit der üblichen lexikographischen Ordnung
zusammenf"allt. Dies erlaubt eine sehr schnelle Kodierung und
Dekodierung von Vektoren in $\overline{L}.$ In \texttt{tn} sind für
eine Liste von Vektoren in $\overline{L}$ die folgenden Methoden für
Listenoperationen implementiert.
\begin{description}
\item[\texttt{insert\_vector (v)}] Füge den Vektor \texttt{v} in die Liste
  ein.
\item[\texttt{delete\_vector (v)}] Lösche den Vektor \texttt{v} aus der
  Liste. 
\item[\texttt{lookup\_vector (v)}] Prüfe, ob der Vektor sich bereits in der
  Liste befindet und liefere dementsprechend \texttt{true} oder \texttt{false}
  zurück. 
\item[\texttt{find\_min}] Liefere den bezüglich der lexikographischen Ordnung
  kleinsten Vektor in der Liste zurück.
\end{description}
Diese Operationen sind aus Effizienzgründen in \texttt{tn}
zum Teil inline implementiert, d.h.\ nicht durch explizite Funktionsaufrufe. Der
Klarheit halber formulieren wir beispielhaft zwei dieser Operationen als
Funktionsaufrufe.  

Die Vektoren werden in einem Array von \texttt{unsigned char}, einem
8-Bit Datentyp kodiert. Dabei geben die führenden $n-3$ Bits die
Position des Bytes im Array an, in welchem der Vektor kodiert wird,
und die letzten drei Bits werden durch die Position in dem
entsprechenden Byte kodiert.

Als erstes Beispiel beschreiben wir die Operation \texttt{insert\_vector(v)}.
 Um das den Vektor \texttt{v}
beschreibende Bit aufzufinden, berechnen wir zun"achst den Offset im
 Array. 
Die führenden $n-3$ Bit von \texttt{v}, interpretiert als
Bin"arzahl, geben die Position des entsprechenden Bytes im Array an
 (im Source Zeile~14
\texttt{l[v>>3]}). In diesem Byte setzen wir das Bit, dessen Position durch
das Bitmuster letzten $3$ Bit von \texttt{v} in der üblichen
Binararithmetik beschrieben wird. Hierzu maskieren wir alle
Bits von \texttt{v} bis auf die letzten $3$ aus (\texttt{(v \& 7)}) und
verwenden das Ergebnis, um das erste Bit um die entsprechende Anzahl zu
verschieben (im Source: \texttt{1 <{}< (v \& 7)}).

\lstset{escapechar=\%}
\begin{lstlisting}[caption={Operationen in $\overline{L}$: \texttt{insert\_vector (v)}}]{}
/* insert vector v in list l */

typedef unsigned int vector_2;
typedef *unsigned char vlist;

/* Insert v in list l. Idea: Lookup the byte 
in vlist l whose offset is binary coded by 
the leading n-3 bits in v.
Now mark in this byte the bit whose number 
is binary coded by the last 3 bits. */

void
insert_vector (vlist l, vector_2 v){
   l[v>>3] |= (1 << (v & 7))\label{insert_vector_fast_decoding}
}
\end{lstlisting}

Bei der folgenden Routine \texttt{find\_min(L)}, die einen Vektor aus
der Liste \texttt{L} 
zurückliefert, müssen wir die Fehlersituation der leeren Liste
abfangen.    Der Nullvektor spielt in
$\overline{L}$ keine Rolle, weil er nicht zur Nachbarbildung verwendet werden
kann. Wir benutzen deshalb sein Kodifikat "`$0$"' als Flag, um
der aufrufenden Funktion mitzuteilen, daß die Liste leer ist.

\lstset{escapechar=\%}
\begin{lstlisting}[caption=%
{Operationen in $\overline{L}$: \texttt{find\_min}}]{}
int dim; /* the dimension of \overline{L} */
typedef unsigned int vector_2;
typedef *unsigned char vlist;
/* Lookup the "smallest" vector in l
   and return it. */

vector_2
find_min (vlist l){
  int i, t;
/* Lookup the first non-empty byte 
   in the array L */
  for (i = 0; i < (1 << (dim - 3)); i++)
    if (l[i] != 0)
      break;

  for (t = 0; t < 8; t++)
    if (((l[i] >> t) & 1) == 1)
      break;
  if (t == 8)
    return 0; /* There is no umarked vector. */                
  else
    return (8 * i + t);
    }
\end{lstlisting}

Um in der Liste den kleinsten Vektor aufzufinden, d.h., in dem Array
das erste gesetzte Bit zu finden, suchen wir die Liste linear durch.
Diese Implementation ist f"ur unsere Zwecke in der Laufzeit
hinreichend effizient, denn wir können jeweils $8$~Bits zugleich
prüfen, in dem wir jeweils ein ganzes Byte betrachten und testen, ob
es gleich $0$ ist.\footnote{Die naheliegende Verbesserung, gleich $32$
  oder mehr Bits durch einen Cast auf einen gr"o"seren Datentyp
  wie Integer zu "uberpr"ufen, f"uhrt nicht zu wesentlichen
  Geschwindigkeitssteigerungen und macht die Portierung auf andere
  Betriebssysteme komplizierter.} Beispielsweise muß in der Dimension
$20$ lediglich ein zusammenh"angender Speicherbereich von höchstens
$128$~kByte durchlaufen werden. Der damit verbundene Aufwand ist
irrelevant gegenüber den anderen Operationen in \texttt{tn} wie
beispielsweise der Basisreduktion. Dieses Durchlaufen des Arrays
geschieht in der folgenden Schleife.
\begin{verbatim}
   for (i = 0; i < (1 << (dim - 3)); i++)
    if (l[i] != 0)
      break;
\end{verbatim}
Nun enth"alt die Variable \texttt{i} den richtigen Offset und wir suchen in dem
\texttt{i}-ten Byte des Arrays das erste gesetzte Bit.
\begin{verbatim}
for (t = 0; t < 8; t++)
  if ((l[i] >> t) & 1)
    break;
\end{verbatim}
Als Ergebnis geben wir den Vektor zurück, dessen Bitmuster der Bin"arzahl $8 i
+ t$ entspricht.

\subsection{Die diskrete Fourier--Transformierte der
  L"angenfunktion}\label{sec:DFT}
\subsubsection{Die diskrete Fourier--Transformation auf endlichen abelschen Gruppen}
\label{sec:DFT_intro}
Das Programm \texttt{tn} enthält eine Vielzahl von Analysefunktionen,
um Invarianten eines Gitters zu bestimmen.
%, siehe Abschnitt~\ref{tn_beschreibung}.
Darunter befindet sich die
bisher wenig untersuchte Funktion  $\mathcal{F}_{L/2L}(\ell),$ die
Fourier-Transformierte 
der L"angenfunktion $\ell.$ Diese Zerlegung der L"angenfunktion, welche
einer Nebenklasse in $L/2L$ die Norm eines  kürzesten in ihr
enthaltenen Vektors zuordnet, in eine
Linearkombination aus irreduziblen Charakteren von $L/2L$ erlaubt bei
vielen Gittern eine recht einfache Beschreibung dieser Funktion. Die
Zerlegung selbst ist abhängig von der Wahl der zugrundeliegenden Basis
in $L$ beziehungsweise in $L/2L;$ die Vielfachheiten der Koeffizienten
in der Linearkombination, das sogenannte Spektrum, sind hingegen unabhängig von
der Wahl einer Basis. Wir beschreiben zunächst die
Fourier-Transformation und ergänzen diese mit einigen Beispielen.

Sei $G$ eine endliche abelsche Gruppe, $G^\vee$ die Menge der
irreduziblen Charaktere von $G$ und seien $f,f'$ komplexwertige Funktionen
auf $G.$ Die irreduziblen Charaktere $G^\vee$ bilden eine
Orthonormalbasis bez"uglich des Skalarproduktes   $ \langle
f,f'\rangle=\frac{1}{|G|}\sum_{g\in G}f(g)\overline{f'(g)}$ im Raum der
Klassenfunktionen $G^\ast$ auf $G.$ Deshalb kann man  $f$ folgendermaßen als
Linearkombination der irreduziblen Charaktere  schreiben, die eine
Orthonormalbasis im Raum der Klassenfunktionen bilden:

\begin{displaymath}
  f=  \sum\limits_{\chi\in G^\vee} \langle \chi,f\rangle\, \chi 
\end{displaymath}
Man nennt $\hat{f}(\chi):= \langle \chi,f\rangle\, =
\frac{1}{|G|}\sum_{g\in G}\chi(g)\overline{f(g)}$ den Fourierkoeffizienten von $\chi$
bez"uglich $f.$

Dadurch wird der Operator $\mathcal{F}$ der
Fourier--Transformation auf der Gruppe $G$ definiert
\begin{displaymath}
  \begin{matrix}
\mathcal{F_G}: & G^\ast &\rightarrow&(G^\ast)^\ast\\
&f&\mapsto&\hat{f}: \chi\mapsto<\chi,f>
  \end{matrix}
\end{displaymath}

Wenn man die Charaktertafel von $G$ als Matrix $(\chi(g))_{\chi\in G^\vee,
  g\in G}$ und $f$ als Spaltenvektor $(f(g))_{g\in G}$ schreibt, dann berechnet
  sich $(\hat{f}(\chi))_\chi$ bis auf einen Faktor als einfache Matrizenmultiplikation durch
  \begin{displaymath}
    (\hat{f}(\chi))_\chi = \frac{1}{|G|}\cdot(\chi(g))_{\chi, g}\, \overline{(f(g))_{g}}\,.
  \end{displaymath}

\subsubsection{Die diskrete Fourier-Transformierte der L"angenfunktion eines
  Gitters}
\label{sec:DFT-L}

Für ein Gitter sei $\ell : L/2L \rightarrow \BR$ mit $\ell(\overline{v})
= \min_{w\in \overline{v}} (w,w)$ als die L"angenfunktion von $L$
definiert.

\paragraph{Der Operator $\mathcal{F}_{L/2L}$}

Die Gruppe $L/2L$ ist eine elementarabelsche Gruppe der Ordnung $2^n$ mit
$n=\rank L.$ Im ersten Schritt verifiziert man, daß sich ihre
Charaktertafel in Form einer Hadamard-Matrix schreiben l"aßt. Diese
Matrix ist bis auf einen skalaren Faktor gleich der Operatormatrix von
$\mathcal{F}_{L/2L}.$ Zun"achst benötigt man eine Notation für die
irreduziblen Charaktere von $L/2L \cong \BF_2^n.$ Die Charaktere
dieser Gruppe lassen sich als mehrfache Tensorprodukte der
Charaktere von $\BF_2$ schreiben.

Die Charaktertafel von $\BF_2$ enth"alt genau zwei Charaktere, den
Hauptcharakter $\chi^+$ und den Signumcharakter $\chi^-.$

\begin{table}[htbp]
\begin{displaymath}
\begin{array}[h]{c|cc}
&0&1\\\hline
\chi^+&1&1\\
\chi^-&1&-1
\end{array}
\end{displaymath}
  \caption{Die Charaktertafel von $\BF_2$}
  \label{tab:chartabf_2}
\end{table}
Die irreduziblen Charaktere von $\BF^n_2$ ergeben sich als $n$-fache
Tensorprodukte der Charaktere $\chi^+$ und $\chi^-.$ Zur Vereinfachung
der Notation wird ein irreduzibler Charakter von $\BF^n$ nur mit einer
Folge von $+$ und $-$ indiziert und nicht als voll ausgeschriebenes
Tensorprodukt, beispielsweise $\chi^{+-++}:= \chi^+\otimes
\chi^-\otimes \chi^+\otimes \chi^+.$ Damit ergeben sich für $\BF_2^1n$
und $\BF_2^2$ und $\BF_2^3$ die Charaktertafeln wie sie in den Tabellen
\ref{tab:chartabf_2}-\ref{tab:chartabf_2_3} beschrieben sind. 

\begin{table}[htbp]
\begin{displaymath}
\begin{array}[h]{c|cccc}
&00&01&10&11\\\hline
\chi^{++}&1&1&1&1\\
\chi^{+-}&1&-1&1&-1\\
\chi^{-+}&1&1&-1&-1\\
\chi^{--}&1&-1&-1&1\\
\end{array}
\end{displaymath}
  \caption{Die Charaktertafel von $\BF_2^2$}
  \label{tab:chartabf_2_2}
\end{table}

\begin{table}[htbp]
\begin{displaymath}
\begin{array}[h]{c|cccccccc}
&000&010&100&110&001&011&101&111\\\hline
\chi^{+++}&1&1&1&1&1&1&1&1\\
\chi^{++-}&1&-1&1&-1&1&-1&1&-1\\
\chi^{+-+}&1&1&-1&-1&1&1&-1&-1\\
\chi^{+--}&1&-1&-1&1&1&-1&-1&1\\
\chi^{-++}&1&1&1&1&-1&-1&-1&-1\\
\chi^{-+-}&1&-1&1&-1&-1&1&-1&1\\
\chi^{--+}&1&1&-1&-1&-1&-1&1&1\\
\chi^{---}&1&-1&-1&1&-1&1&1&-1\\
\end{array}
\end{displaymath}
  \caption{Die Charaktertafel von $\BF_2^3$}
  \label{tab:chartabf_2_3}
\end{table}

Das Programm \texttt{tn} berechnet die Fourier-Transformierte der
L"angenfunktion $\ell$. Bei der Berechnung einer Vielzahl von Gittern stellt
sich heraus, da"s sich die Koeffizienten in der Zerlegung in irreduzible
Charaktere auf nur wenige Werte bei gro"sen Vielfachheiten konzentrieren.

Die folgenden Beispiele sollen diese Beobachtung illustrieren und zu neuen
Fragestellungen Anla"s geben. 

\subsubsection{Beispiele}
\label{sec:dft_beispiel}
\begin{enumerate}
\item Sei $L=(\BZ^2,q(.,.))$ mit
\begin{math}
G(q)=
\begin{pmatrix}
  1&0\\
  0&2
\end{pmatrix}
\end{math}
Dann ist
\begin{eqnarray*}
  (\hat{\ell}(g))_g = \frac{1}{|L/2L|} \cdot(\chi(g))_{\chi, g}\,
{(\ell(g))_{g}} = \\
\frac{1}{4} 
  \begin{pmatrix}
1&1&1&1\\
1&-1&1&-1\\
1&1&-1&-1\\
1&-1&-1&1\\
\end{pmatrix}
\cdot
\begin{pmatrix}
  0\\
  2\\
  1\\
  3\\
\end{pmatrix}
= \frac{1}{4}
\begin{pmatrix}
  6\\
  -4\\
  -2\\
  0
\end{pmatrix}
.
\end{eqnarray*}
Also l"a"st sich ${\ell}$ in die folgende Linearkombination irreduzibler
Charaktere zerlegen
\begin{displaymath}
  \ell = \frac{1}{4}(6\chi^{++}-4\chi^{+-}-2\chi^{-+}+0\chi^{--} ).
\end{displaymath}

\item Sei $A_2=(\BZ^2,q(.,.))$ mit
\begin{math}
G(q)=
\begin{pmatrix}
  2&-1\\
  -1&2
\end{pmatrix}.
\end{math}
Dann zerf"allt $\ell$ analog des vorhergehenden Beispiels in die 
Linearkombination  
\begin{displaymath}
\ell=\frac{1}{4}(6\chi^{++} -2\chi^{-+} -2\chi^{+-} -2\chi^{--}).
\end{displaymath}
\item F"ur gr"o"sere $n$ zeigt sich bei der Zerlegung von $\ell$ f"ur die
  Gitter $A_n$ eine Verteilung der Koeffizienten auf nur
  wenige Werte.
  \begin{enumerate}
  \item Koeffizientenverteilung f"ur $A_{12}$
    \begin{center}
    \begin{tabular}[t]{c|c}
      Koeffizient&Vielfachheit\\\hline
      26624&1\\
      -2048&13\\
      0&4082\\
    \end{tabular}
  \end{center}
%$f(A_{12}) = 13 \times -2048 + 4082 \times 0 + 1 \times 26624 $
  \item Koeffizientenverteilung f"ur $A_{13}$

    \begin{center}
    \begin{tabular}[t]{c|c}
      Koeffizient&Vielfachheit\\\hline
      57344&1\\
      -4096&14\\
      0&8177\\
    \end{tabular}
  \end{center}
%   $f(A_{13}) = 14 \times -4096 + 8177 \times 0 + 1 \times 57344 $
  \item  Koeffizientenverteilung f"ur $A_{14}$

    \begin{center}
    \begin{tabular}[t]{c|c}
      Koeffizient&Vielfachheit\\\hline
      122880&1\\
      -8192&15\\
      0&16358\\
    \end{tabular}
  \end{center}
%$f(A_{14}) = 15 \times -8192 + 16368 \times 0 + 1 \times 122880 $
  \end{enumerate}
\end{enumerate}

\section{$2$-Nachbarn bei gerader Determinante}
\label{sec:tn_even_det}
Bisher haben wir die $2$-Nachbarbildung nur bei geraden Gittern mit
ungerader Determinante untersucht. Den Grund für diese Einschr"ankung
finden wir im Beweis von
Lemma~\ref{lem:construction_of_neighbours}. $L_v$ mu"s bei der
Nachbarbildung ein echtes Untergitter von $L$ sein; dies ist im
allgemeinen aber nicht gegeben, wie wir im Beispiel $L=2M$ für ein
ganzzahliges Gitter $M$ sehen. Dann gilt für jeden Vektor $v\in L,$
daß $L=L_v$ ist. F"ur viele Vektoren $v$ ist $L\neq L_v$ aber auch
bei gerader Determinante erf"ullt und damit  ist die Konstruktion
eines Nachbarn ebenfalls m"oglich.

Für die Konstruktion des $2$-Nachbarschaftsgraphen eines Gitters
$L$ mit gerader Determinante, das aus den Vektoren $b_1,\cdots,b_n$
erzeugt wird, geht man im wesentlichen
genau so vor wie im ungeraden Fall. Man berechnet die Menge aller Orbits von
$\overline{O(L)}$ auf $\overline{L}$ und w"ahlt aus jedem Orbit einen
Repr"asentanten $\overline{v}$ aus. Dann sucht man einen Basisvektor $b_i$ von
$L$ mit $(\overline{v},\overline{b_i})=1$ und $\overline{v}\neq
\overline{b_i}.$ Diese Eigenschaften sind unabh"angig von der Auswahl des
Repr"asentanten $v.$ Wie in Lemma~\ref{lem:construction_of_neighbours}~i)
gezeigt, ist bei ungerader Determinante die Existenz eines solchen $b_i$
gesichert, im geraden Fall jedoch nicht. 
Findet man keinen Basisvektor $b_i$
mit diesen Eigenschaften, dann verwirft man diesen Orbit und geht zum
n"achsten über. 

Dieses Vorgehen garantiert natürlich nicht eine erfolgreiche Klassifikation
eines ganzen Geschlechtes, führt in der Regel aber in größeren Dimensionen
(etwa ab Dimension $5$ bis $8$) zum Ziel. In fast allen geraden Gittern
$L$ höherer Dimension finden sich gen"ugend Orbits von Vektoren $v\in L$ mit
$L_v\neq L,$ welche die Existenz von $2$-Nachbarn oft
sicherstellen. Wir verfolgen 
hier mit \texttt{tn} einen pragmatischen Ansatz. \texttt{tn} berechnet den
Nachbarschaftsgraphen, und eine erfolgreiche Verifikation mit der Maßformel
sichert die Vollst"andigkeit des Ergebnisses
(s.~Abschnitt~\ref{sec:verifikation}). Wie oben bereits beschrieben, sind
jedoch Beispiele von Gittern ohne $2$-Nachbarn auch in hohen Dimensionen
konstruierbar.

Als prominentes Beispiel f"ur eine erfolgreiche Klassifikation in einem
Geschlecht mit gerader Determinante führen wir in
Abschnitt~\ref{sec:barnes_wall} die 
Klassifikation des Geschlechtes des Barnes-Wall-Gitters an.

\section{Ergebnisse}
\label{sec:tn_results}
Die wichtigsten Ergebnisse, die mit Hilfe von \texttt{tn} erzielt wurden, sind
in dem Artikel \cite{scharlau98:_class} dokumentiert. Wir geben hier eine
kurze Zusammenfassung und einige neuere Resultate an.

\subsection{$\ell$-modulare Gitter}
\label{sec:modular_lattices}
Ein $n$-dimensionales ganzzahliges Gitter $L$ heißt \emph{modular vom Level
  $\ell,$} wenn es ein $\ell\in\BN$ gibt mit $L\cong \sqrt{\ell} L^\#$
\cite{quebbemann95}. Es folgt unmittelbar, daß $\det L = \ell^{n/2};$
insbesondere ist $n$ gerade, wenn $\ell$ quadratfrei ist. Viele der
  aus der Literatur bekannten
Gitter sind modular, so z.B. das $12$-dimensionale Coxeter-Todd-Gitter,
das $16$-dimensionale Barnes-Wall-Gitter und die geraden unimodularen
Gitter $E_8,$ das Gosset-Gitter, und  $\Lambda_{24},$  das
  Leech-Gitter \cite{bible}. 
 Aus der Theorie der
Modulformen ergibt sich für modulare Gitter eine obere Schranke für das
Minimum. Gitter, die dieses theoretische Maximum annehmen, heißen
\emph{extremal.} Für eine Einordnung und Diskussion dieses Themas
  verweisen wir auf 
  den Übersichtsartikel
\cite{scharlau98:_extrem_lattic}. Die Existenz und die Eindeutigkeit modularer
Gitter ist eine zentrale Frage in dieser Theorie.  Einige der extremalen
Gitter realisieren die besten bekannten Packungsdichten oder würden
sie im Falle der Existenz übertreffen.
  
Für eine Berechnung des Nachbarschaftsgraphen mit \texttt{tn}
benötigen wir ein Startgitter aus dem zu bestimmenden Geschlecht.
Außer bei $\ell=5$ haben wir bei den unten aufgeführten F"allen (hier
gilt jeweils $\ell
\equiv 3 \pmod{4}$) mit dem direkten Produkt aus $n/2$ Gittern
 in der Dimension $2$ mit Grammatrix
\begin{math}
  \begin{pmatrix}
    \ell/2&1\\1&2
  \end{pmatrix}
\end{math}
begonnen.

Wir zitieren hier eine Liste mit den aus der Theorie der Modulformen
berechneten möglichen maximalen Minima.
\begin{table}[htbp]
  \begin{displaymath}
    \begin{array}[b]{ll|rrrrrrrr}
      & n   & 4 & 6 & 8 & 10 & 12 & 14 & 16 & 18 \\
      \hline
      \ell=3:&\min & 2 & 2 & 2 & 2 & 4 & 4 & 4 & 4 \\
      \hline
      \ell=5:&\min & 2 &   & 4 &   & 4 &   & 6 &   \\
      \hline
      \ell=7:&\min & 2 & 4 & 4 & 4 & {6} & 6 & 6 & 8 \\
      \hline
      \ell=11:&\min & 4 & 4 & 6 & 6 & {8} & 8 &  &   \\
    \end{array}
  \end{displaymath}
  \caption{Maximal mögliche Minima von modularen Gittern (aus \cite{scharlau98:_class})}
  \label{tab:minima_modularer_gitter}
\end{table}
\subsubsection{Der Level $\ell=3$}
Die Klassenzahlen für die Dimensionen $2,4,6,8,10,12$ sind jeweils
$1,1,1,2,3,10;$ alle diese Gitter sind modular. Bis auf das
Coxeter-Todd-Gitter in Dimension $12$ mit Minimum $4$ sind alle
anderen Gitter reflektiv, d.h., sie enthalten ein Wurzelsystem von
vollem Rang (\cite{scharlau96:_reflec}). Eine Methode, das Geschlecht
des Coxeter-Todd-Gitters über die Konstruktion der möglichen
Wurzelsysteme zu bestimmen, und die daraus folgende Klassifikation
wurden bereits in \cite{scharlau95:_coxet_todd_lattic,scharlau00:_coxet_todd_lattic} ohne Detailrechnungen angegeben.

Mit dem Programm \texttt{tn} l"aßt sich aufgrund des
Konstruktionsverfahrens nicht nur eine Auflistung aller Gitter eines
Geschlechts berechnen, sondern auch der vollst"andige
Nachbarschaftsgraph beschreiben. Die folgende Tabelle enth"alt die
Inzidenzmatrix des Nachbarschaftsgraphen des Geschlechtes des
Coxeter-Todd-Gitters, bei welchem die Kanten mit der Anzahl der
verschiedenen Orbits gewichtet wurden, die unter der jeweiligen
Operation von $\overline{O(L)}$ auf $L/2L$ zu einer
Nachbarschaftsbildung führen.

\begin{table}[htbp]
\begin{center}
  \begin{math}
\begin{array}[h]{r|rrrrrrrrrr}
  &  K_1&  K_2&  K_3&  K_4&  K_5&  K_6&  K_7&  K_8&  K_9& K_{10}\\\hline
K_1 &1&1&&&&&&&&\\
K_2&1&5&1&1&1&1&&&&\\
K_3&&1&4&1&&&1&&&\\
K_4&&1&1&7&1&&1&1&1&\\
K_5&&1&&1&5&1&1&&1&\\
K_6&1&&&1&1&1&&&&\\
K_7&&&1&1&1&&6&1&&1\\
K_8&&&&1&&&1&4&1&\\
K_9&&&&1&1&&&1&4&1\\
K_{10}&&&&&&&1&&1&2\\
\end{array}
\end{math}
\caption{Der Nachbarschaftsgraph des Geschlechtes des
    Coxeter-Todd-Gitters mit den Bezeichnungen aus \cite{scharlau95:_coxet_todd_lattic}}
  \label{tab:graph_coxeter_todd}
\end{center}
\end{table}

In der Dimension $14$ gibt es $29$ Gitter in
diesem Geschlecht, darunter ein extremales mit Minimum $4.$ Seine
Automorphismengruppe von der Ordnung $2^7\cdot 3^6\cdot 7\cdot 13$ ist (die
modulo der Spiegelung durch den Ursprung einfache Gruppe) $G_2(3)\cdot 2$
(\cite{scharlau98:_class}, Proposition~3.2). Eine explizite
Konstruktion dieses Gitters 
wird im $\mathbb{ATLAS}$ \cite{atlas}, S.~60 beschrieben.

\subsubsection{Der Level $\ell=5$}
Solche Gitter existieren nur in durch $4$ teilbaren Dimensionen. In
der Dimension $12$ existiert ein interessantes Geschlecht mit
Determinante $5^6,$ welches aus $48$ Gittern besteht, darunter $40$ modularen
und $4$ extremalen.  Die Klassifikation dieses Geschlechtes
wurde unabh"angig von G.~Nebe (unveröffentlicht) erzielt.

\subsubsection{Der Level $\ell=7$}
Das Hauptergebnis des Klassifikationsprojektes in
\cite{scharlau98:_class} ist der Beweis der Nichtexistenz eines
vermuteten Gitters.
In der
 Dimension $12$ existiert kein extremales Gitter mit der Determinante $7^6.$
Dies ist das kleinste bekannte Beispiel für die Nichtexistenz eines möglichen
extremalen Gitters. 
% Es ist im Gegensatz zum Geschlecht $11^6$
% (siehe den folgenden 
% Abschnitt~\ref{11h6}) kein Beweis ohne Computerunterstützung bekannt.

\subsubsection{Der Level $\ell=11$}
In den Dimensionen $4,6,8$ existieren jeweils $3,5,31$ Gitter, alle sind
modular und in jedem Geschlecht gibt es ein eindeutiges extremales Gitter. 

In der Dimension $10$ gibt es neben dem Craig-Gitter ein weiteres extremales
(Beobachtung von H.-G.~Quebbemann, \cite{quebbemann95}) alle weiteren der
insgesamt $297$ Gitter besitzen ein kleineres Minimum. Siehe
Abschnitt~\ref{11h6} für den $12$-dimensionalen Fall.

\subsection{Das Geschlecht des Barnes-Wall-Gitters
  $\Lambda_{16}$}\label{sec:barnes_wall} 

Das Gitter $\Lambda_{16}$ wurde 1959 von Barnes und Wall  konstruiert
\cite{barnes59:_some_abelian}  und ist mit
Minimum $4$ und Determinante $256$ das dichteste bekannte Gitter in
Dimension $16.$ Die Bezeichnung $\Lambda_{16}$ ist durch seine Konstruktion
als geschichtetes Gitter motiviert (siehe \cite{bible}, Chapt.~6).
 Das Geschlecht des Barnes-Wall-Gitters
wurde  von R.~Scharlau und B.B.~Venkov  \cite{MR95e:11073}
klassifiziert. Analog zur Klassifikation der geraden, unimodularen Gitter in
Dimension $24$ von B.B.~Venkov in \cite{venkov:_even_unimod_dimen_lattic}
werden die Gitter im Geschlecht des Barnes-Wall-Gitters über ihr Wurzelsystem
bestimmt.

Obwohl es sich hier um ein Geschlecht von gerader Determinante handelt, setzen
wir \texttt{tn} mit den in Abschnitt~\ref{sec:tn_even_det} beschriebenen
Modifikationen ein. Ausgehend von einer Erzeugermatrix für $\Lambda_{16}$
untersuchen wir zun"achst beispielhaft
die Existenz eines Nachbarvektors und pr"asentieren 
anschließend die Klassifikation und den $2$-Nachbarschaftsgraphen des
kompletten  Geschlechtes des Barnes-Wall-Gitters.

Eine Erzeugermatrix $A$ für $\Lambda_{16}$ entnehmen wir
\cite{bible}, S.~130.
\begin{small}
  \begin{displaymath}
    A=\frac{1}{\sqrt{2}}\left(
      \begin{array}{cccccccccccccccc}
        4& & & & & & & & & & & & & & & \\
        2&2& & & & & & & & & & & & & & \\
        2&0&2& & & & & & & & & & & & & \\
        2&0&0&2& & & & & & & & & & & & \\
        2&0&0&0&2& & & & & & & & & & & \\
        2&0&0&0&0&2& & & & & & & & & & \\
        2&0&0&0&0&0&2& & & & & & & & & \\
        2&0&0&0&0&0&0&2& & & & & & & & \\
        2&0&0&0&0&0&0&0&2& & & & & & & \\
        2&0&0&0&0&0&0&0&0&2& & & & & & \\
        2&0&0&0&0&0&0&0&0&0&2& & & & & \\
        1&1&1&1&0&1&0&1&1&0&0&1& & & & \\
        0&1&1&1&1&0&1&0&1&1&0&0&1& & & \\
        0&0&1&1&1&1&0&1&0&1&1&0&0&1& & \\
        0&0&0&1&1&1&1&0&1&0&1&1&0&0&1& \\
        1&1&1&1&1&1&1&1&1&1&1&1&1&1&1&1
      \end{array}
    \right)
  \end{displaymath}
\end{small}
Wir betrachten den elften und zwölften Zeilenvektor aus der Matrix $A$ und
setzen $v=\frac{1}{\sqrt{2} }(1,1,1,1,0,1,0,1,1,0,0,1,0,0,0,0).$ $v$ ist ein
Minimalvektor, also primitiv und insbesondere nicht in $2\Lambda_{16}$
enthalten. Das Skalarprodukt von $v$ mit dem elften Basisvektor
$w=\frac{1}{\sqrt{2}}(2,0,0,0,0,0,0,0,0,0,2,0,0,0,0,0)$ ist gleich $1.$ Also
ist $w\notin (\Lambda_{16})_v,$ und damit ist gezeigt, daß
$(\Lambda_{16})_v\neq \Lambda_{16}.$ Da $(v+2w,v+2w)=4+4+16=24 \equiv 0
\pmod 8,$ erhalten wir nach Lemma~\ref{lem:integral_even} mit
$\Lambda_{16}(v+2w)$ einen geraden Nachbarn von $\Lambda_{16}.$

\subsubsection{Der Nachbarschaftsgraph}
Wir geben in Tabelle \ref{tab:grap_BW_gitter} den Nachbarschaftgraphen
des Geschlechtes des Barnes-Wall-Gitters an. Die Bezeichnung der
Gitter folgt der Klassifikation in \cite{MR95e:11073}.  Die
Vielfachheiten in der Inzidenzmatrix geben die Anzahl der Orbits an,
die zu einer Nachbarbildung führen.

\begin{table}[htbp]
\begin{center}
\begin{sideways}
\begin{small}
\begin{tabular}[h]{r|rrrrrrrrrrrrrrrrrrrrrrrrr}
  &  1&  2&  3&  4&  5&  6&  7&  8&  9& 10& 11& 12& 13& 14& 15& 16& 17& 18& 19& 20& 21& 22& 23& 24\\\hline
1&        1&1&&&&&&&&&&        &&&&&&&&&&&&\\
2&        1&7&1&2&2&&1&&1&&&&&        &&&&&&&&&&\\
3&&1&5&&2&1&&1&&1&&&&&&        &&&&&&&&\\
4&&2&&3&1&&&1&1&&1&&&&&&&        &&&&&&\\
5&&2&2&1&19&2&2&1&3&4&2&2&1&&&1&&&        &&&&&\\
6&&&1&&2&8&1&&&2&&1&1&1&1&&&&&&        &&&\\
7&&1&&&2&1&5&&1&&1&1&1&&&1&&&&1&&&        &\\
8&&&1&1&1&&&6&&1&2&&&&1&&1&1&&&&&&        \\
9&&1&&1&3&&1&&4&1&1&1&&&&1&&&&&&&&        \\
10&&&1&&4&2&&1&1&17&1&3&2&1&1&1&&1&1&&&&&        \\
11&&&&1&2&&1&2&1&1&10&1&1&&1&1&1&3&&&1&&&        \\
12&&&&&2&1&1&&1&3&1&8&&1&1&1&&&1&1&&&&        \\
13&&&&&1&1&1&&&2&1&&10&1&&3&&1&1&1&1&1&&        \\
14&&&&&&1&&&&1&&1&1&4&&&&&&1&&&1&        \\
15&&&&&&1&&1&&1&1&1&&&7&&&2&&&1&&1&        \\
16&&&&&1&&1&&1&1&1&1&3&&&5&&1&1&1&&1&&        \\
17&&&&&&&&1&&&1&&&&&&1&1&&&&&&        \\
18&&&&&&&&1&&1&3&&1&&2&1&1&8&1&&2&&1&        \\
19&&&&&&&&&&1&&1&1&&&1&&1&3&&&1&&        \\
20&&&&&&&1&&&&&1&1&1&&1&&&&5&1&1&&1        \\
21&&&&&&&&&&&1&&1&&1&&&2&&1&3&1&1&        \\
22&&&&&&&&&&&&&1&&&1&&&1&1&1&2&&1        \\
23&&&&&&&&&&&&&&1&1&&&1&&&1&&2&        \\
24&&&&&&&&&&&&&&&&&&&&1&&1&&1
\end{tabular}
\end{small}
\end{sideways}
  \caption{Der Nachbarschaftsgraph des Geschlechtes des Barnes-Wall-Gitters}
  \label{tab:grap_BW_gitter}
\end{center}
\end{table}
\subsection{Neuere Ergebnisse anderer Autoren}
\label{sec:weiterentwicklungen}

\subsubsection{$11$-modulare Gitter in Dimension $12$}\label{11h6}
Ein extremales $11$-modulares Gitter mit Minimum $8$ in Dimension $12$
w"are geringfügig dichter als das dichteste bekannte Gitter in dieser
Dimension, das Coxeter-Todd-Gitter. Das Maß dieses Geschlechtes liegt
bei ungef"ahr $15096.99$ (\cite{scharlau98:_class}, S.~748), aber die
meisten Gitter besitzen relativ große Automorphismengruppen, so daß
man aufgrund der Maßformel (siehe Gleichung~\ref{eq:massformel} auf
Seite~\pageref{eq:massformel})) von vielen hunderttausend oder
Millionen Gittern in diesem Geschlecht ausgehen kann.  Eine
vollst"andige Klassifikation mit \texttt{tn} erscheint weder sinnvoll
noch aussichtsreich.

In \cite{nebe96:_nonex} zeigen G.~Nebe und B.B.~Venkov, da"s kein solches
extremales Gitter existiert. In dem Beweis, der sich im wesentlichen auf
analytische Methoden und Siegel-Modulformen vom Grad~$2$ stützt, wird
\texttt{tn} als Hilfsmittel benutzt, um ein Erzeugendensystem des
Vektorraums der vorkommenden 
Thetareihen zu konstruieren.

\subsubsection{Die Nachbarmethode im hermiteschen Fall}

In \cite{schiemann98:_class} stellt A.~Schiemann das Programm \texttt{hn} vor,
da"s die Nachbarmethode im hermiteschen Fall anwendet.  
  Die erfolgreichen
Klassifikationen sind in \cite{schiemann98:_class} beschrieben.

\subsubsection{Extremale Gitter}

In dem Übersichtsartikel  über extremale
Gitter von R.~Scharlau und R.~Schulze-Pillot
\cite{scharlau98:_extrem_lattic} wurden die Ergebnisse mit
\texttt{tn} und dem oben erw"ahnten \texttt{hn} berechnet.

\subsubsection{Diverses}

In der Doktorarbeit \cite{mischler96:_un_bz_f} werden von M.~Mischler
einige Geschlechter mit \texttt{tn} bestimmt.

\section{Verifikation der Ergebnisse}
\label{sec:verifikation}

Bei allen Klassifikationsproblemen stellt sich die Frage nach der Korrektheit
und der Vollst"andigkeit des Ergebnisses.

Das sogenannte Maß eines Geschlechts  $m(\mathcal{G})$ kann einerseits
mit analytischen Methoden hergeleitet werden
(s.~\cite{conway88:_low}), genügt andererseits nach Siegels Maßformel 
der Identit"at
\begin{equation}
  m(\mathcal{G}) = \sum_{[L]\in\mathcal{G}}
  \frac{1}{|O(L)|}\,,\label{eq:massformel}
  \end{equation}
wobei die Summe über die Menge aller mit \texttt{tn} in $\mathcal{G}$ gefundenen Isomorphieklassen $[L]$ gebildet wird.

Diese Testfunktion erlaubt analog einer Quersumme
in den oben zitierten Beispielen eine sehr
verl"aßliche Überprüfung der berechneten Klassifikationen. 
Die Ma"sformel ist sensitiv  gegenüber
fehlenden  Isomorphieklassen, einer
fehlerhaften Berechnung der (Ordnung der) Automorphismengruppen und
inkorrekten 
Isometrietests. Eine fehlerhafte Berechnung eines Geschlechtes, die von der
Maßformel nicht aufgedeckt wird, ist nur möglich bei gleichzeitigem Auftreten
von zwei Fehlern, die sich bei der Berechnung der Maßformel gegenseitig
aufheben. Die an diesem Test beteiligten Programmteile sind die Generierung
der Nachbarn, die Berechnung ihrer Automorphismengruppen und Isometrietests
gegen die bisher gefundenen Gitter. Arbeiten zwei dieser drei Teile korrekt,
so wird ein Fehler im dritten durch die Maßformel erkannt. 
Die Algorithmik des Isometrietests und der Berechnung der
Automorphismengruppe sind sich in großen Programmteilen ähnlich, weil
der Isometrietest mathematisch vergleichbar mit dem Auffinden eines
nichttrivialen Automorphismus ist.  In diesen Programmteilen befinden
sich aber eine ganze
Anzahl von Plausibilit"atstests, die Inkonsistenzen entdecken k"onnen.

Die wenigen bereits bekannten Klassifikationen (s.o.) konnten sowohl für
ungerade als auch für gerade Determinanten erfolgreich reproduziert werden.
Die prominente Klassifikation der Niemeier-Gitter
\cite{niemeier73:_defin_formen_dimen_diskr} konnte mit \texttt{tn} nicht
durchgeführt werden. Der Grund liegt in der algorithmischen Komplexit"at der
Reduktion von Gitterbasen und in dem Umstand, daß die Automorphismengruppen in
\texttt{autom} als Permutationsgruppen auf der Menge der kurzen Vektoren
realisiert sind. Die Basis eines Nachbargitters, das mit dem oben
beschriebenen Algorithmus konstruiert wird, enth"alt zumeist einige l"angere
Vektoren. In der Routine \texttt{autom} müssen alle Vektoren des Gitters bis
zu der L"ange des l"angsten Basisvektors berechnet werden. In der Dimension $24$
ist eine vollst"andige Aufz"ahlung der kurzen Vektoren praktisch nur für
Vektoren der L"ange $2$ und $4$ möglich. Der LLL-Algorithmus reduziert aber
nicht alle erzeugten Gitterbasen auf eine Basis, die nur Vektoren der L"ange
kleiner oder gleich $4$ enthalten. So bleibt die Klassifikation von
Niemeier ein 
Meilenstein in der Theorie der ganzzahligen quadratischen Formen, für die es
auch heute kein algorithmisches "Aquivalent gibt.

\chapter{"Uber die Zerlegung von Gittern}\label{cha:mmlll}
\section{"Uberblick}

Im Mittelpunkt dieses Kapitels steht ein  Algorithmus nach
einer geometrischen Idee M.~Knesers zur
Zerlegung eines Gitters 
in eine orthogonale Summe von unzerlegbaren Untergittern. Eine
vereinfachte Variante dieser Methode ermöglicht eine effiziente 
Konstruktion einer 
Gitterbasis aus einem großen Er\-zeu\-gen\-den\-sys\-tem.  Eine frühere Version
dieses Kapitels ist  
 in einen gemeinsamen Artikel \cite{hemkemeier98} mit Frank
Vallentin eingeflossen.

In diesem Kapitel ist es zweckm"aßig, ein Gitter $L$ als Untermodul im
$n$-di\-men\-si\-ona\-len euklidischen Raum $E$ versehen mit dem
Standardskalarprodukt $(.,.)$ zu betrachten und den Aspekt der
quadratischen Form als beschreibende Invariante eines Gitters in den
Hintergrund treten zu lassen. Das Gitter $L$ soll immer von vollem Rang sein,
d.h., $\dim \Rspan{L}=n.$  Im folgenden bezeichnen wir  
hier mit der Norm eines Vektors $x\in E$ seine euklidische Norm
$||x||=\sqrt{(x,x)}.$

\subsection{Das Zerlegungsproblem}
\label{sec:zerlegungsproblem}
Bei vielen algorithmischen Konstruktionen von Gittern hat man nicht
gen"ugend Kontrolle "uber die Eigenschaften und Invarianten des zu
konstruierenden 
Gitters. Insbesondere ist es f"ur Klassifikationszwecke interessant,
ob ein Gitter die orthogonale Summe zweier
niederdimensionalerer Gitter ist. In dieser Situation ben"otigen wir
einen effizienten  Algorithmus, der ein Gitter in eine orthogonale
Summe von Teilgittern zerlegt.

\begin{definition}\label{def:indecomposable_lattice}
  Ein nichttriviales Gitter $L$ hei"st \emph{zerlegbar,} falls es
  (echte) Untergitter $L',L''\subset L$ gibt mit $L=L'+ L''$ und
  $(L_1,L_2) = 0,$ andernfalls \emph{unzerlegbar.} 
\end{definition}
\begin{remark}
  Wir schreiben $L=L'\oplus L'',$ falls $L$ in $L'$ und $L''$
  zerlegbar ist. In diesem Kapitel werden wir das Symbol $\oplus$
  unter Strapazierung der Notation \textbf{aus\-schliess\-lich} f"ur eine
  innere direkte, orthogonale Summe verwenden.
\end{remark}

Die Existenz einer Zerlegung in unzerlegbare Untergitter
 folgt unmittelbar durch Induktion
über den Rang von $L.$ In einem nicht konstruktiven Beweis zeigte Eichler
1952 in \cite{eichler52} die Existenz einer bis auf die Reihenfolge
der Summanden eindeutigen Zerlegung in unzerlegbare Untergitter. 
M.~Kneser zeigte kurz darauf in einem elegantem Beweis \cite{kneser54}, da"s
diese Zerlegung effektiv berechenbar ist. Mittels eines Erathostenessiebes
werden Gittervektoren eliminiert, die nicht in einem unzerlegbaren
Untergitter liegen. Damit 
erh"alt man ein Erzeugendensystem f"ur die unzerlegbaren Untergitter. In
Abschnitt~\ref{sec:knesers_method} geben wir eine Version von Knesers Beweis
f"ur Erzeugendensysteme an, welche die benutzten geometrischen Methoden
illustriert. 

Bei der Umsetzung von Knesers Beweis in einen Algorithmus stößt man 
auf zwei Probleme. Das Aussieben der zerlegbaren Vektoren ist ein 
aufwendiger Vorgang, dessen Laufzeit im wesentlichen  quadratisch
von der 
Gr"o"se des Erzeugendensystems abh"angt. Des weiteren bildet die
unzerlegbaren Vektoren im 
allgemeinen selbst ein sehr gro"ses Erzeugendensystem f"ur die unzerlegbaren
Untergitter. Die Berechnung einer Gitterbasis aus einem Erzeugendensystem ist
zwar ein algorithmisch gut verstandenes Problem, in der Praxis jedoch f"ur
gro"se Erzeugendensysteme  aufwendig. Die bekannten Algorithmen wie
der MLLL und die Bestimmung der Hermite-Normalform sind für
Erzeugendensysteme konstruiert, die nicht wesentlich größer als der
Rang des erzeugten Gitters sind. In der oben
beschriebenen Situation führt dies in der Regel zu einem relativ schlechten
Laufzeitverhalten.

In Abschnitt~\ref{sec:improved_algorithm} wird ein  schnellerer
Algorithmus f"ur das Zerlegungsproblem als die Implementation
von Knesers Methode angegeben; insbesondere erfolgt das Aussieben
der erzeugenden Vektoren in linearer Laufzeit. 

\subsection{Minimale Erzeugendensysteme}
\label{sec:size_of_generating_systems}
Das zentrale Ergebnis in diesem
Kapitel ist ein Überdeckungssatz für Er\-zeu\-gen\-den\-sys\-teme, der
eine obere Schranke  f"ur die Gr"o"se eines minimalen
Erzeugendensystems aus Vektoren mit einer beschr"ankten L"ange angibt. 

In Gittern mit nicht zu kleinem Minimum $M$ (betrachtet im Verh"altnis
zur Determinante) besteht jedes minimale Erzeugendensystem nur aus
wenig mehr als $n$ Vektoren.  Der
"Uberdeckungssatz~\ref{satz:covering_theorem} quantifiziert dieses
Ergebnis genauer: Ist $S\subset L$ ein minimales Erzeugendensystem
von $L$ mit 
Vektoren, die nicht l"anger als $B\in \BR$ sind, dann ist
$|S|\leqslant  n +
\log_2(n!(\frac{B}{M})^n).$ 

Die algorithmische Hürde dieser Aufgabe rührt vom Fehlen eines
Basiserg"anzungssatzes und eines Steinitz'schen Basisaustauschsatzes für
die Modulsituation her. Man kann sich dies leicht an einem Erzeugendensystem
veranschaulichen, in dem kein Vektor im Erzeugnis der anderen liegt,
siehe etwa das Beispiel~\ref{bsp:all_vectors_needed} auf
Seite~\pageref{bsp:all_vectors_needed}.

\subsection{Der MLLL bei gro"sen Erzeugendensystemen}
\label{sec:MLLL_large_generating_systems}
Die in Kapitel~\ref{sec:improved_algorithm} entwickelten Ideen für
eine effiziente Lösung des Zerlegungsproblems
werden
in Kapitel~\ref{sec:mmlll} eingesetzt, um einen effizienten
Algorithmus zu konstruieren, der aus einem gro"sen Erzeugendensystem
eine Gitterbasis berechnet. Die Methode besteht aus einem
Meta-Algorithmus, der als Kern bekannte Routinen wie den MLLL oder die
Berechnung der Hermite-Normalform benutzt. Diese Kernalgorithmen
werden wiederholt aufgerufen, jedoch ausschließlich für
minimale, also sehr kleine Erzeugendensysteme. Für den
Meta-Algorithmus können bei großen Erzeugendensystemen sch"arfere
asymptotische Schranken gezeigt werden als für die zugrundeliegenden
Kernalgorithmen. Die praktische Relevanz erkannt man an den
Laufzeitgewinnen, die in den graphischen Auswertungen zweier Experimente
(s.~S.~\pageref{fig:dim20_100}, Abb.~\ref{fig:dim20_100} und
Abb.~\ref{fig:dim20_1000}) veranschaulicht werden.

\section{Unzerlegbare Vektoren}
\label{sec:unzerlegbare_vektoren}
Bevor wir uns der konkreten Berechnung einer Zerlegung von $L$ in
unzerlegbare Untergitter zuwenden, betrachten wir qualitative
Kriterien für die Zerlegbarkeit bzw.~{}
Unzerlegbarkeit von Vektoren. F"ur ein Gitter mit der
Zerlegung $L=L_1 \oplus \ldots\oplus L_r$ sollen unzerlegbare Vektoren
zwei Eigenschaften gen"ugen:
\begin{enumerate}
\item Wenn $x\in L_i,\, y\in L_j, i\neq j$ unzerlegbare Vektoren sind,
  dann sei $x+y\in L_i\oplus L_j$ zerlegbar.
\item Die Menge der unzerlegbaren Vektoren in $L_i, 1\leqslant i
  \leqslant r$ erzeugt $L_i.$
\end{enumerate}

Der Begriff der Unzerlegbarkeit von Vektoren
wurde zuerst von M.~Kneser definiert \cite{kneser54}
und sp"ater von O'Meara \cite{omeara80} verfeinert.\footnote{Wir danken
  M.~Kneser für seinen freundlichen Hinweis auf die Arbeit von O'Meara.}

\begin{definition}\label{def:indecomposable_vector}
  Sei $v\in L.$
  \begin{enumerate}
  \item $v$ hei"st \emph{orthogonal zerlegbar,} falls Vektoren $x,y\in
    L\setminus 0$ existieren mit $v=x+y$ und $(x,y)=0,$ andernfalls
    \emph{orthogonal unzerlegbar.}
  \item $v$ hei"st \emph{spitz zerlegbar,} falls Vektoren $x,y\in
    L\setminus 0$ existieren mit $v=x+y$ und $(x,y)\geqslant 0,$
    andernfalls \emph{spitz unzerlegbar.}
  \item $v$ hei"st \emph{linear zerlegbar} oder einfach
    \emph{zerlegbar,} falls Vektoren $x,y\in L$ existieren mit $v=x+y$
    und $||v|| > ||x||\geqslant ||y||,$ andernfalls \emph{(linear)
      unzerlegbar.}
  \end{enumerate}
\end{definition}

\begin{lem}[\cite{omeara80}]
  Sei $v\in L,$ dann gilt:
  \\$v$ unzerlegbar $\Rightarrow v$ spitz unzerlegbar $\Rightarrow v$
  orthogonal unzerlegbar.  
\end{lem}
\begin{proof}
  Sei $v$ orthogonal zerlegbar, dann ist $v$ offenbar auch spitz
  zerlegbar. Sei $v$ nun spitz zerlegbar, d.h., es gibt $x,y\in
  L\setminus 0$ mit $v=x+y$ und $(x,y)\geqslant 0.$ Dann folgt aus
  $||v||= \sqrt{x^2 +2(x,y) + y^2}> %\geqslant
  ||x||,$ da"s $v$ zerlegbar ist.
\end{proof}

O'Meara  gibt eine Charakterisierung
orthogonal und spitz zerlegbarer Vektoren an.
\begin{satz}[\cite{omeara80} 3.6, 3.7]\label{satz:spitz_unzerlegbar} Sei $v\in L.$
  \begin{enumerate}
  \item $v$ ist genau dann orthogonal zerlegbar, wenn  ein
    $w\in L$ existiert mit $w\neq \pm v,\, v-w\in 2L$ und $||v||=
    ||w||.$ 
  \item\label{spitz}  $v$ ist genau dann spitz zerlegbar, wenn  ein
    $w\in L$ existiert mit $w\neq \pm v,\, v-w\in 2L$ und
    $||v||\geqslant     ||w||.$ 
  \end{enumerate}
\end{satz}
\begin{proof}
  Wir folgen O'Mearas Beweis und zeigen zun"achst Aussage~ii). 
  Sei $v\in L.$
  \begin{align*}
    &v \text{ spitz zerlegbar}\\
    \Longleftrightarrow\quad& \exists y\in L\setminus \{0,v\}\quad
    (y,v-y)\geqslant 0\\& (\text{Nun sei } w:=v-2y \text{ bzw.\ }
    y:=\frac{1}{2}(v-w)).\\ 
    \Longleftrightarrow\quad& \exists w\in L\setminus \{\pm v\}\quad
    (\frac{1}{2}(v-w),\frac{1}{2}(v+w)) \geqslant 0 \text{ und }
    v-w\in 2L\\
    \Longleftrightarrow\quad& \exists w\in L\setminus \{\pm v\}\quad
    (v,v)-(w,w) \geqslant 0 \text{ und }
    v-w\in 2L\\
    \Longleftrightarrow\quad& \exists w\in L\setminus \{\pm v\}\quad
    ||v||\geqslant ||w|| \text{ und }
    v-w\in 2L\\
  \end{align*}
  
  Aussage~i) erhalten wir, wenn wir im vorangehenden Beweis "`spitz"'
  durch "`orthogonal"' und alle "`$\geqslant$"' durch "`$=$"'
  ersetzen.
\end{proof}
\begin{figure}[htbp]
  \begin{center}
  \includegraphics[%width=0.8\textwidth,
  scale=0.4,bbllx=40mm,bburx=200mm,bblly=80mm,bbury=200mm,clip=]{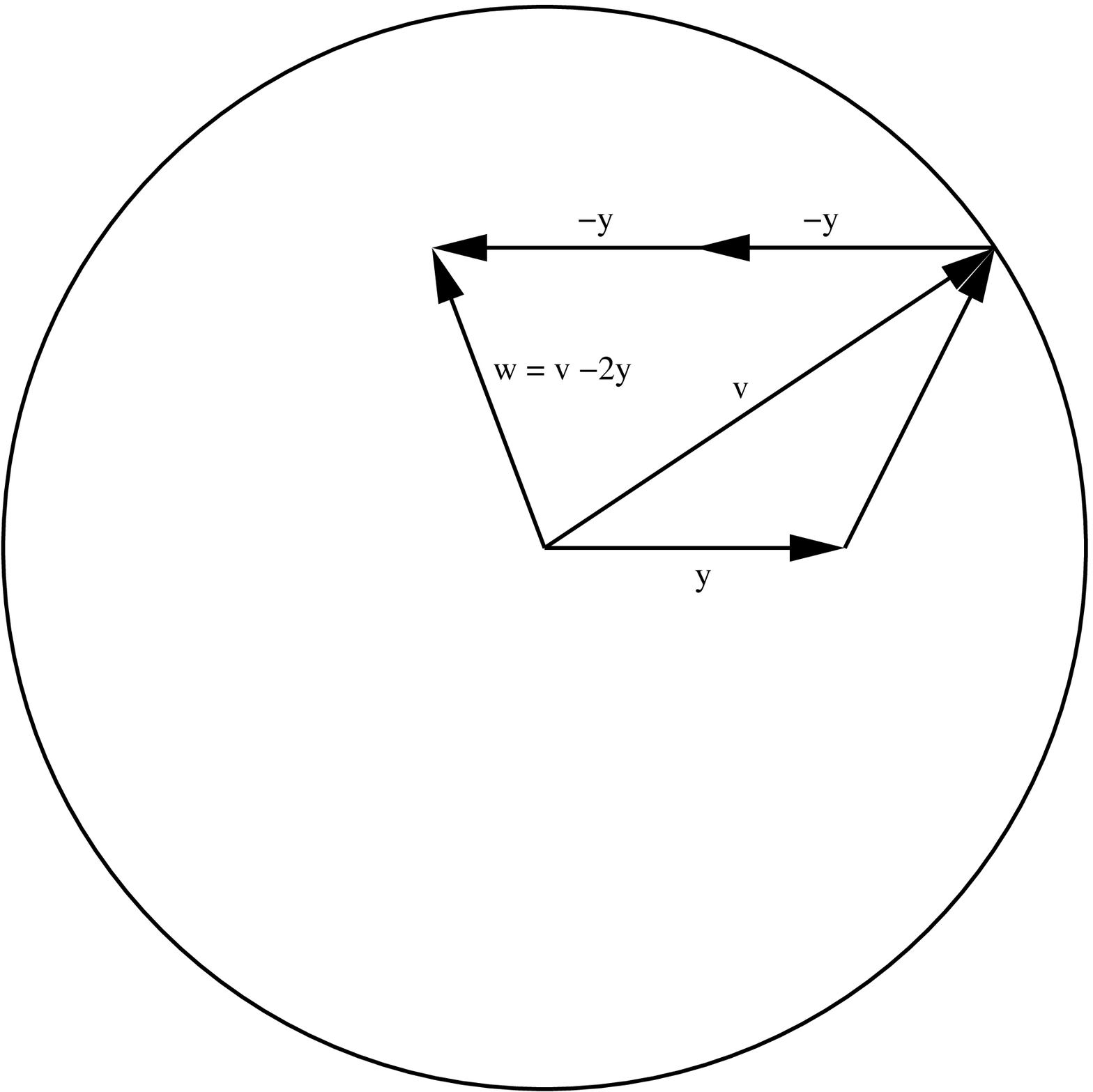}\\
  \small F"ur $v,y\in L$ mit $(y,v-y)\geqslant 0$ ex.\ ein $w\in L$
  mit $||w||\leqslant ||v||$ und $v-w\in 2L.$
  \caption{Illustration zum Beweis von  Satz~\ref{satz:spitz_unzerlegbar} ii) "`$\Rightarrow$"'}
  \label{fig:indec_2}
\end{center}
\end{figure}

Damit hat O'Meara eine bemerkenswerte Eigenschaft von vorono"irelevanten
Vektoren bewiesen, allerdings ohne den Zusammenhang zur
Vorono"itheorie zu erw"ahnen.\footnote{Siehe die Definition
  \ref{def:voronoirelevant} auf 
  Seite~\pageref{def:voronoirelevant}.} Unseres Wissens ist dieser
simple Zusammenhang bis heute nicht formuliert
worden.
\begin{kor}\label{kor:voronoirelevant}
  Ein Vektor $0\neq x\in L$ ist genau dann spitz unzerlegbar,
  wenn er vorono"irelevant ist. 
\end{kor}
\begin{proof}
  Nach Satz~\ref{satz:spitz_unzerlegbar} ist $x$ genau dann spitz
  unzerlegbar, wenn $\pm x$ die einzigen k"urzesten Vektoren in der
  Nebenklasse $x+2L$ sind. Die Behauptung folgt nun unmittelbar aus
  einem Kriterium von Vorono"i, hier~Satz~\ref{voronoi_relevante_vektoren}, S.~\pageref{voronoi_relevante_vektoren}.  
\end{proof}

Aus Satz~\ref{satz:spitz_unzerlegbar} folgt ebenfalls, da"s es
h"ochstens $2\cdot(2^n-1)$ spitz unzerlegbare (und damit linear
zerlegbare) Vektoren gibt (\cite{omeara80} 3.17, 3.18). Eine einfache
geometrische Überlegung zeigt bereits, da"s lange Vektoren immer zerlegbar
sind. 

\begin{lem}\label{lem:covering_decomp}
  Sei $R$ der "Uberdeckungsradius\footnote{Der "Uberdeckungsradius ist der
    maximal m"ogliche Abstand eines Punktes in E von einem Gitterpunkt in L,
    siehe Definition~\ref{def:covering_radius},
    S.~\pageref{def:covering_radius}.} von $L$ und $v\in L$ mit $||v||> 2R.$
  Dann ist $v$ zerlegbar.
\end{lem}
\begin{proof}[Beweis I:]
  Sei  $w\in L$ mit $||w-\frac{1}{2}v||\leqslant R.$ Dann gilt
  \begin{displaymath}
     ||w||\leqslant ||w-\frac{1}{2}v|| + ||\frac{1}{2}v|| \leqslant
     R+\frac{1}{2}||v|| < ||v|| 
  \end{displaymath}
  und
  \begin{displaymath}
     ||w-v||\leqslant ||w-\frac{1}{2}v||+||-\frac{1}{2}v||\leqslant
     R+\frac{1}{2}||v|| < ||v|| .
  \end{displaymath}
  Also ist $v$ zerlegbar, weil $v$ die Summe der beiden echt k"urzeren
  Vektoren $w$ und $v-w$ ist.
\end{proof}
\begin{figure}[htbp]
  \begin{center}
    \includegraphics[]{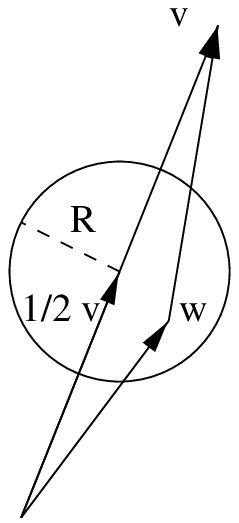}
    \caption{Lemma \ref{lem:covering_decomp}: Lange Vektoren sind zerlegbar}
    \label{fig:lange_vektoren_sind_zerlegbar}
  \end{center}
\end{figure}

Mit einem Argument aus der Vorono"itheorie folgt die Behauptung unmittelbar.
\begin{proof}[Beweis II:]
(F. Vallentin)
  Ist $||v||> 2R,$ so folgt aus $\frac{1}{2}v\notin V_0(L),$ daß $v$ nicht
  vorono"irelevant ist und mit Satz~\ref{satz:spitz_unzerlegbar}
  zerlegbar sein 
  muß.
\end{proof}

\begin{remark}
  Der Begriff der orthogonalen Zerlegbarkeit von Vektoren wurde von
  M.~Kneser analog zum Begriff der Zerlegbarkeit von Gittern definiert
  \cite{kneser54}.  In generischen Gittern sind alle Vektoren
  orthogonal unzerlegbar. Ein einfaches Beispiel hierf"ur ist das von
  den Vektoren $(1,0), (\theta,1)$ mit $\theta^2\notin \mathbb{Q}$
  erzeugte Gitter.  Insbesondere sind in einem eindimensionalen Gitter
  alle nichttrivialen Vektoren orthogonal unzerlegbar. O'Mearas Definition der
  linearen Zerlegbarkeit von Vektoren erscheint uns intuitiver.
\end{remark}

\section{Knesers Methode}
\label{sec:knesers_method}
Sei $B\in \BR,$ so da"s $S:=\{v_1,\cdots,v_s\} = \{v\in L\setminus 0
\mid ||v|| \leqslant B\}$ ein Erzeugendensystem von $L$ ist. Auf der
Menge der unzerlegbaren Vektoren $I$ in $S$ definieren wir den
Orthogonalit"atsgraph $\Gamma=(I,E)$ mit $E=\{\,\{v,w\}\in I\times
I\mid (v,w)\neq 0\}.$ Wir zerlegen $\Gamma$ in seine
Zusammenhangskomponenten $\Gamma_i=(I_i,E_i), 1\leqslant i \leqslant
r$ und setzen $L_i=\Zspan{I_i}.$
\begin{satz}[\cite{kneser54}]
  Die Gitter $L_i, 1\leqslant i \leqslant r$ sind unzerlegbar und die
  orthogonale Zerlegung 
  \begin{equation}
    L=L_1\oplus \ldots\oplus L_r\label{eq:mmlll_decomposition}
  \end{equation}
  ist eindeutig bis auf die Reihenfolge der Summanden.
\end{satz}
\begin{proof}
  Sei $v\in S\setminus \Zspan{I}$ und $||v||$ minimal mit dieser
  Eigenschaft. Dann ist $v$ wegen $v\notin I$ 
zerlegbar und es gibt $x,y\in S$ von echt
  k"urzerer Norm als $v$  mit $v=x+y.$ Mindestens einer dieser beiden
  Vektoren ist nicht in $\Zspan{I},$ sonst w"are auch
  $v\in\Zspan{I}.$ Dies widerspricht jedoch der Minimalit"at von
  $||v||,$ somit ist $S\subseteq \Zspan{I}$ gezeigt. Daraus folgt
  die erste Behauptung $L=\sum_{i=1}^r L_i.$
  
  $(I_i,I_j)=0$ f"ur $i\neq j$ impliziert $(L_i,L_j)=0.$ Daraus folgt
  mit der Definitheit von $(.,.)$ und der vorhergehenden Aussage, 
da"s $$L=L_1\oplus \ldots\oplus L_r.$$
  
  Seien $L',L''\subseteq L_i$ f"ur ein beliebiges, aber festes $i$
mit $L_i=L'\oplus L''.$ F"ur $v\in I_i$
  gilt entweder $v\in L'$ oder $v\in L'',$ ansonsten w"are $v$
  zerlegbar.  Ohne Beschr"ankung der Allgemeinheit 
 sei $v\in L'.$ Wegen $(L',L'')=0$ sind $L'\cap
  \Gamma$ und $L''\cap \Gamma$ zwei nicht zusammenh"angende Teilmengen 
 von $\Gamma.$ Also ist $L''=0$ und $L'=L_i.$ Damit ist
  gezeigt, da"s $L_i$ unzerlegbar ist. 

Da jedes $L_i$ bereits durch
  einen in ihm enthaltenen unzerlegbaren Vektor eindeutig bestimmt
  ist, folgt die Eindeutigkeit der Zerlegung bis auf Reihenfolge der
  Summanden. 
\end{proof}

Wir beschreiben nun Knesers Methode zur Konstruktion von $\Gamma:$
\begin{description}
  
\item[Initialisierung] Sei $\{b_1,\cdots,b_n\}$ eine Basis von $L.$
  Setze $B=\max_i ||b_i||$ und berechne $S=\{v\in L\setminus 0\mid
  ||v||\leqslant B\}.$
\item[Konstruiere $I$] F"ur jedes Paar $u,v\in S$ mit   $u+v\in S$
  markiere $u+v$ als zerlegbar. Sei $I$ die Menge aller nicht
  als zerlegbar
  markierten  Vektoren.
\item[Konstruiere $E$] F"uge $\{u,v\}$ f"ur alle $u,v\in I$ in
  $E$ ein, falls $(u,v)\neq 0.$
\end{description}

Bei diesem Verfahren m"ussen wir jedes Paar $\{u,v\}$ f"ur alle
$u,v\in S$ untersuchen. Diese $\binom{s}{2}\in O(s^2)$ Operationen
sind f"ur gro"se $s$ sehr aufwendig.

\section{Der inkrementelle Algorithmus}
\label{sec:improved_algorithm}
Um für ein $L_i$ in Gleichung \ref{eq:mmlll_decomposition} eine Basis
effizient zu
bestimmen, ist es nicht notwendig, alle unzerlegbaren Vektoren in
$L_i$ 
zu berechnen. Es gen"ugt, ein Erzeugendensystem in den $L_i$
aufzufinden, welches einerseits aus relativ wenigen Vektoren besteht
und andererseits erlaubt, schnell zu verifizieren, da"s sich alle
unzerlegbaren Vektoren in den daraus erzeugten Gittern befinden. Dies motiviert den
folgenden Algorithmus~\ref{algo:improved_algorithm}.

\begin{algorithm}[h]
\caption{Zerlegung des Gitters $L$}
\label{algo:improved_algorithm}
\begin{algorithmic}[1]
\STATE \emph{Input:} $B\in\BR$ mit einem Erzeugendensystem $S=\{v\in
L\setminus 0\mid ||v||\leqslant B\}$ von  $L.$ 

\STATE \emph{Output:} Unzerlegbare Untergitter $L_i$ mit
$L=L_1\oplus\ldots\oplus L_r.$

  \STATE {\sl // 1.\ Initialisierung.}
  \STATE W"ahle $v \in S$ mit  $||v||$ minimal.
  \STATE $S \leftarrow S\backslash\{v\}$, $k \leftarrow 1$, $L_k  \leftarrow
  {\Zspan{v}}$

  \STATE {\sl // 2.\ F"uge sukzessiv Vektoren aus $S$ in die $L_i$ ein.}
  \WHILE {$S \neq \emptyset$} 
  \STATE  W"ahle $v \in S$ mit minimalem $||v||$.
  \STATE $S \leftarrow   S\backslash\{v\}$.
  \IF {$v \not\in L_1 \oplus\ldots\oplus L_k$}
  \STATE {\sl // $v$ is ,,erg"anzend'', insbesondere unzerlegbar.}
  \STATE $J \leftarrow \{ i \in   \{1, \cdots, k\} \mid \pi_i(v) \neq
  0\}$
  \STATE $\tilde{L} \leftarrow \BZ v   + \bigoplus_{i\in J} L_i$ 
  \STATE \textsl{// Ordne die Liste der $L_i$ neu.} 
  \STATE $\{L_1,\cdots,L_{k-|J|}\} \leftarrow \{L_j
  \mid j\notin J\}, \quad L_{k-|J|+1} \leftarrow \tilde{L}, \quad k\leftarrow
  {k-|J|+1}$ 
  \ENDIF
  \ENDWHILE
\end{algorithmic}
\end{algorithm}

Die Korrektheit dieses Algorithmus zeigt 
\begin{satz}\label{theorem:correctness_of_algorith_1}
  Sei $L$ ein $n$-dimensionales Gitter auf $E,$  $B\in\BR,$ so da"s
  die Menge
  $S=\{v_1,\cdots,v_s\} := \{v\in L \mid ||v|| \leqslant
  B\}$ ein Erzeugendensystem von $L$ enth"alt. Dann berechnet
  Algorithmus~\ref{algo:improved_algorithm} die Zerlegung von $L$ in
  unzerlegbare Untergitter $L_i,$ d.~h., 
  $L=L_1\oplus\ldots\oplus L_r.$
\end{satz}
F"ur ein gegebenes System von paarweise orthogonalen Untergittern $L_1
\oplus\cdots\oplus L_k$ sei $\pi_i$ die orthogonale Projektion auf den
$\BR$-Vektorraum, der von $L_i$ aufgespannt wird. $I$ und $\Gamma$
seien wie oben definiert.

\begin{proof}
  Sei $L_1\oplus\ldots\oplus L_k$ mit den folgenden Eigenschaften
  bereits berechnet:
  \begin{enumerate}
  \item $L=\Zspan{S} + \sum_{i=1}^k L_i,$
  \item $L_i$ ist unzerlegbar f"ur $1\leqslant i\leqslant k,$
  \item $(L_i,L_j)=0$ f"ur $i\neq j.$
  \end{enumerate}
  Diese Eigenschaften sind nach der Initialisierung in Zeile~5 von
  Algorithmus~\ref{algo:improved_algorithm} erf"ullt.

  Wir zeigen nun die Schleifeninvarianz (Zeile 7--16) dieser Eigenschaften.
  Sei $v\in S$ ein Vektor mit minimaler Norm in $S$ wie in Zeile~8 gew"ahlt.
Ist
  $v\in L_1\oplus\ldots\oplus L_k,$ so k"onnen wir $v$ verwerfen und die
  Eigenschaften i) -- iii) bleiben erhalten.  Sei nun $v\notin
  L_1\oplus\ldots\oplus L_k.$ Einen solchen Vektor nennen wir
  ,,erg"anzend''\footnote{Die Begriffsbildung ist lediglich lokal,
  d.h.\ nur im jeweiligen Schleifendurchlauf
  sinnvoll, weil diese 
  Eigenschaft von der Menge der bereits abgearbeiteten Vektoren in $S$
abh"angt und
  nicht eine Invariante des Vektors $v$ ist.}.
  \begin{enumerate}
    \renewcommand{\labelenumi}{Zu \roman{enumi}):}
  \item $v$ wird hinzugef"ugt und es gilt 
    \begin{displaymath}
      L=\Zspan{S\setminus \{v\}} + (\BZ v + \sum_{i=1}^k L_i).
    \end{displaymath}
    Insbesondere gilt $L=L_1\oplus\ldots\oplus L_r$ nachdem alle
    Vektoren in $S$ abgearbeitet sind.
  \item $v$ ist unzerlegbar, sonst g"abe es $x,y\in S$ von echt
    k"urzerer Norm als $v$ mit $v=x+y.$ Nach Voraussetzung w"aren
    $x,y$ jedoch schon in einer vorherigen Schleife abgearbeitet
    worden, woraus $x+y= v\in L_1\oplus\ldots\oplus L_k$ folgen w"urde.
    (Insbesondere folgt aus der orthogonalen Summe $v=(v-\pi_i(v)) +
    \pi_i(v), 1\leqslant i\leqslant k,$ da"s es keine nichttriviale
    Projektion von $v$ auf $\Qspan{L_i}$ gibt, welche in $L_i$ liegt.)
    Wir setzen $J = \{ j \in \{1, \cdots, k\} \mid \pi_j(v) \neq 0\}$
    und w"ahlen Vektoren $v_j\in L_j, j\in J$ mit $(v_j,v)\neq 0.$
    Dann ist der Teilgraph $\{(v_j,v)\}_{j\in J}\subseteq
    \Gamma$ ein Stern und damit zusammenh"angend, also ist
    $\bigoplus_{j\in J} L_j + \BZ v$ unzerlegbar. 
  \item $\bigoplus_{i\in I\setminus J} L_i \oplus (\bigoplus_{j\in J}
    L_j + \BZ v)$ ist eine orthogonale Zerlegung, weil $\pi_i(v) = 0$
    f"ur $i\notin J.$
  \end{enumerate}
  
  Damit ist die Korrektheit von
  Algorithmus~\ref{algo:improved_algorithm} gezeigt.
\end{proof}

\subsection{Datenstrukturen und Laufzeitanalyse}
\label{sec:improved_algo_data_struct}
In diesem Abschnitt beschreiben wir f"ur die Realisierung von
Algorithmus~\ref{algo:improved_algorithm} geeignete Datenstrukturen.
Wir werden in den S"atzen \ref{theorem:runtime_of_algorith_1} und
\ref{theorem:mmlll} Laufzeitanalysen f"ur die hier entwickelten
Algorithmen aufgrund dieser Datenstrukturen vornehmen.  Im folgenden
fassen wir $E$ als Vektorraum aus Spaltenvektoren auf.

Eine zentrale Rolle in diesem Algorithmus spielen die orthogonalen
Zerlegungen $L_1 \oplus\ldots \oplus L_k,$ die in jedem
Schleifendurchlauf konstruiert werden.
Jedes Teilgitter $L_i$ wird
durch eine Gitterbasis $v_{i1},\cdots,v_{ir_i}\subset E$
repr"asentiert. Des weiteren ben"otigen wir eine (m"oglicherweise
leere) Vektorraumbasis $w_1,\cdots,w_{n-l}$ des orthogonalen Komplements
$\Rspan{L_1,\cdots,L_k}^\perp,\, l = n- \sum_{i=1}^k r_i.$ Nun setzen wir
\begin{displaymath}
  A= 
  \begin{pmatrix}
    v_{11}&\cdots&v_{kr_k}&w_1&\cdots&w_{n-l}
  \end{pmatrix}^{-1}\in \RMat{n}.
\end{displaymath}

F"ur jedes $x\in E$ gilt $v\in L_1 \oplus\ldots \oplus L_k$ genau
dann, wenn die ersten $n-l$ Koeffizienten von $Av$ in $\BZ$ liegen und
die restlichen $l$ gleich $0$ sind. Wir k"onnen somit in Zeile~$10$
 des Algorithmus
durch eine
Matrix-Vektor-Multiplikation testen, ob der Vektor $x$ in dem
bereits konstruierten Untergitter liegt. Ist das der Fall, so
"uberspringen wir ihn und gehen zum n"achsten
Schleifendurchlauf. Andernfalls berechnen wir wie folgt
das neue Gitter $\tilde{L}$ mit
einer seiner Basen, ordnen die Liste unserer unzerlegbaren Teilgitter
$L_i$ neu und berechnen die Indexmenge $I$ und dann
 die Matrix $A$ in der
folgenden Weise neu: $A$\label{mlll:neuberechnung_von_A} ist eine 
Diagonalmatrix mit den Bl"ocken 
\begin{math}
  \begin{pmatrix}
    v_{i1}&\cdots&v_{ir_i}
  \end{pmatrix}^{-1}
\end{math}
auf der Hauptdiagonalen. F"ur jeden der ersten $l$ Koeffizienten von
$Av,$ welcher nicht in $\BZ$ liegt, suchen wir den entsprechenden
Block auf (etwa $j$) und f"ugen $j$ in $J$ ein. Anschlie"send
konstruieren wir eine Basis von $\tilde{L}=\sum_{j\in J} L_j$ mit
einer Variante\footnote{Beispielsweise mit einer der in
  \cite{pohst87:_lll} oder \cite{buchmann89:_comput} beschriebenen
  Modifikationen. Siehe Abschnitt~\ref{sec:mmlll} f"ur eine
  ausf"uhrliche Diskussion der Alternativen.} des bekannten
LLL-Algorithmus \cite{lenstra82:_factor}. Die Berechnung einer Basis
von $\tilde{L}$ kann effizient erfolgen, weil $\tilde{L}$ ein
Erzeugendensystem $\{v_{ji}\mid j\in J,1\leqslant i\leqslant
r_j \}\cup \{v\}$ aus h"ochstens $\sum_{j\in J} r_j \leqslant n+1$
Vektoren besitzt.  Nun ordnen wir die Liste der Gitter $\tilde{L},L_i,
i\notin J$ neu und  
indizieren sie von $1$ bis $k-|J|+1.$ Schlie"slich berechnen wir $A$
mittels Gau"stransformationen neu.

Die Effizienz dieses Algorithmus beruht auf der Tatsache, da"s sich
die meisten $v\in S$ in einer bereits  konstruierten Zerlegung $L_1
\oplus\ldots \oplus L_k$ befinden und die aufwendige Berechnung des
Gitters $M$ und seiner Basis im Verh"altnis zur Mächtigkeit von  $S$ nur
selten notwendig ist. Im folgenden werden wir diese Behauptung genauer
quantifizieren. Dazu benutzen wir haupts"achlich Minkowskis zweiten
Fundamentalsatz.

\subsubsection{Ketten von Untergittern}

\begin{satz}[\cite{minkowski96:_geomet_zahlen}]\label{minkowskis_second_theorem}
  Sei $L$ ein Gitter auf dem $n$-dimensionalen euklidischen Raum $E$
  mit den sukzessiven Minima $\lambda_1(L), \cdots, \lambda_n(L)$
  und $B_n=\{x\in E\mid ||x||\leqslant 1\}$ der $n$-dimensionale Ball
  mit Radius $1.$ Dann gilt
  \begin{equation}
    \frac{2^n}{n!} \sqrt{\det L} \leqslant
    \lambda_1(L)\cdot\ldots\cdot\lambda_n(L) \cdot \vol B_n\leqslant
    2^n\sqrt{\det L}.\label{eq:minkowski}
  \end{equation}
\end{satz}

Im Zentrum dieses Abschnittes steht der folgende
\begin{satz}\label{proposition:bounded_chain}
  Sei $L$ ein Gitter mit Minimum $M$ auf dem $n$-dimensionalen
  euklidischen Raum $V.$ Das Untergitter $L'\subseteq L$ sei ebenso wie
  $L$ von 
Vektoren erzeugt, die alle nicht l"anger als
  $B\in \BR$ sind. Weiter gebe es eine endliche Kette von Untergittern
  $L_i\subseteq L$ mit
  \begin{equation}
    L'=L_1\subset \cdots \subset L_t=
    L\label{eq:lattice_chain}
  \end{equation}
  Dann gilt f"ur die L"ange $t$ der Kette:
  \begin{enumerate}
  \item Ist $\rank L'=\rank  L,$ dann ist $$ t\leqslant
    \log_2(n!\left(\frac{B}{M}\right)^n).$$ 
  \item Wenn alle $L_i$ aus Vektoren nicht l"anger als $B$ erzeugt werden,
    dann  ist $$ t\leqslant n+
    \log_2(n!\left(\frac{B}{M}\right)^n).$$ 
  \end{enumerate}
\end{satz}
\begin{proof}
  Zu i): Sei $s_i=[L:L_i]$  der Index von $L_i$ in $L, 1\leqslant
    i\leqslant t.$ Dann gilt nach der Determinanten-Index-Formel
    $s_i^2=[L:L_i]^2 = \frac{\det L_i}{\det L},$ da"s
  \begin{displaymath}
    \sqrt{\det L'/\det L} = [L:L'] = s_1 >\cdots> s_t
    = [L:L] = 1.
  \end{displaymath}
  $s_{i+1}$ ist wegen $s_{i+1}\cdot [L_{i+1}:L_i] = s_i$ ein Teiler
  von $s_i.$ Also ist die maximale Anzahl von Gittern in der Kette
  (\ref{eq:lattice_chain}) durch die l"angstm"ogliche Teilerkette von
  $\sqrt{\det L'/\det L}$ nach oben beschr"ankt. Der kleinstm"ogliche
  Teiler ist jeweils $2,$ also erhalten wir f"ur $t$ eine obere Schranke von
  $\log_2 \sqrt{\det L'/\det L}.$

  Nach Satz \ref{minkowskis_second_theorem} gilt f"ur $L'$
  \begin{displaymath}
    \sqrt{\det L'} \leqslant \frac{n!}{2^n}\vol B_n\cdot
    \lambda_1(L')\cdot\ldots\cdot\lambda_n(L')
    \leqslant \frac{n!}{2^n}\vol B_n B^n
  \end{displaymath}
  und für $L$
  \begin{displaymath}
    \sqrt{\det L} \geqslant \frac{1}{2^n}\vol B_n\cdot
    \lambda_1(L)\cdot\ldots\cdot\lambda_n(L)
    \geqslant \frac{1}{2^n}\vol B_n M^n.
  \end{displaymath}
  Daraus folgt
  \begin{displaymath}
    s_1=[L:L']=\sqrt{\frac{\det L'}{\det L}}
    \leqslant n!\left(\frac{B}{M}\right)^n.
  \end{displaymath}

  Zu ii) Wir nehmen ohne Beschr"ankung der Allgemeinheit an, da"s $\rank
  L_{i+1} - \rank L_{i} \leqslant 1$ f"ur alle $1\leqslant i < n$ ist.
  Ansonsten verfeinern wir die Kette und f"ugen neue Gitter hinzu, bis sich
  zwei in der Kette benachbarte Gitter jeweils im Rang um h"ochstens $1$
  unterscheiden. Dann zeigen wir die Behauptung f"ur die verfeinerte,
  l"angere Kette. 

  Nun werden die Gitter umnumeriert, um  sie mit ihrem Rang zu
  indizieren. Die Kette bestehe nun aus Gittern 
  $L_{ij}\subseteq L,\, 1\leqslant j\leqslant m_i$ vom Rank
  $i,$ erzeugt von Vektoren nicht l"anger als $B\in \BR,$ so da"s
  \begin{equation}
    L_{01}\subset L_{11}\subset L_{12}\subset\cdots \subset L_{n,m_n-1}\subset
  L_{nm_n} = L.\label{eq:complete_chain}%\label{eq:lattice_chain} 
  \end{equation}

  F"ur $1\leqslant i\leqslant n$ sei $z_i\in L_{i1}\setminus
  \Rspan{L_{i-1,m_{i-1}}}$ mit $||z_i||\leqslant B.$ Ein solches $z_i$
  existiert, weil alle Gitter $L_{ij}$ aus Vektoren nicht l"anger als $B$
  erzeugt werden.  Wir setzen $F=\langle z_1,\cdots,z_n\rangle_\BZ.$ Die Menge
  $\{ z_1,\cdots,z_n\}$ ist nach Konstruktion $\BR-$linear unabh"angig, also
  ist $\rank F=n.$ Nun addieren wir $F$ zu jedem Gitter in der
  Kette~(\ref{eq:complete_chain}) hinzu. Damit erhalten wir die folgende Kette
  von Gittern, die alle den Rang $n$ besitzen.
  \begin{equation}
    \label{eq:blow_up_chain}
    (F=)L_{01}+F\subseteq\ldots \subseteq L_{ij}+F\subseteq\ldots
\subseteq L_{nm_n}+F(=L).
  \end{equation}
  Man verifiziert leicht, da"s 
 $L_{ij} +F\neq L_{i,j+1}+F,\, j < r_i$ f"ur alle $1\leqslant
 i\leqslant n$ gilt. Dann kann in der Kette (\ref{eq:blow_up_chain})
 Gleichheit nur an den Positionen $L_{i-1,r_{i-1}}+F\subseteq L_{i1}+F, 1\leq
 i\leq n,$ also höchstens $n$-mal gelten. 
Wir wenden die in Teil i) gezeigte Aussage  auf
die Kette der verschiedenen Gitter in (\ref{eq:blow_up_chain}) an und
 haben damit gezeigt, da"s die Kette
 (\ref{eq:complete_chain}) und damit auch die Kette in Behauptung ii) nicht mehr als $n+\log_2(n!(\frac{B}{M})^n)$
 Gitter enth"alt. 
\end{proof}

\subsubsection{Ein "Uberdeckungssatz f"ur Erzeugendensysteme}
Wir fassen diese Ergebnisse in einem "`"Uberdeckungssatz"' zusammen,
der die Existenz von kleinen Erzeugendensystemen garantiert.
\begin{definition}
  Ein Erzeugendensystem $S$ eines Gitters heißt \emph{minimal,}
  wenn f"ur alle $S'$ gilt
  \begin{displaymath}
    S'\subset S \Rightarrow \Zspan{S'}\subset\Zspan{S}.
  \end{displaymath}
\end{definition}
\begin{satz}\label{satz:covering_theorem}
  Sei $L$ ein Gitter mit Minimum $M$ auf dem $n$-dimensionalen euklidischen
  Raum $E.$ Sei $S$ ein minimales 
  Erzeugendensystem von $L$ mit Vektoren, die nicht l"anger als $B\in
  \BR$ sind. Dann  ist $$|S|\leqslant  n +
  \log_2(n!(\frac{B}{M})^n).$$
\end{satz}
\begin{proof}
  Sei $S=\{s_1,\cdots,s_t\}.$ Wegen der Minimalit"at von $S$ gilt
  $\Zspan{s_1,\cdots,s_i}\neq   \Zspan{s_1,\cdots,s_{i+1}}$ f"ur alle
  $1\leqslant i<t.$ Die Behauptung folgt aus der Anwendung von
  Satz~\ref{proposition:bounded_chain}~ii) auf die Kette
  \begin{displaymath}
    \Zspan{s_1}\subset \Zspan{s_1,\cdots,s_i} \subset
    \Zspan{s_1,\cdots,s_{t}}=\Zspan{S}. 
  \end{displaymath}
\end{proof}
\subsubsection{Laufzeitanalyse von Algorithmus
  \protect\ref{algo:improved_algorithm} }
Eine  asymptotische Schranke f"ur die Anzahl der "`teuren"'
Updateoperationen in Algorithmus~\ref{algo:improved_algorithm}, bei
denen die Basis von $\tilde{L}$ und die Matrix $A,$ wie auf
Seite~\pageref{mlll:neuberechnung_von_A} beschrieben, neu berechnet 
werden müssen, erhalten wir in 
\begin{kor}
    \label{kor:number_of_update_steps}
    Die Anzahl der Updateoperationen in
    Algorithmus~\ref{algo:improved_algorithm} für die Zerlegung eines
    Gitters $L$ ist in
    $O(n\log\frac{nB}{M}).$ 
\end{kor}
\begin{proof}
  In Algorithmus~\ref{algo:improved_algorithm} werden sukzessive
  Zerlegungen von Teilgittern $L_1\oplus \ldots\oplus L_k\subseteq L$
  konstruiert. Diese Teilgitter erf"ullen die Voraussetzungen von
Satz \ref{proposition:bounded_chain}~ii). Damit ist die Anzahl
  der erg"anzenden Vektoren, die eine Updateoperation notwendig
  machen, nicht gr"o"ser als 
  $n+\log_2(n!(\frac{B}{M})^n).$ Auf diese Schranke wenden wir die
  Stirling'sche Formel $n!\approx \sqrt{2\pi n}(\frac{n}{e})^n$ an und
  erhalten eine obere Absch"atzung, die in 
 $O(n\log\frac{nB}{M})$ liegt.
\end{proof}

Ein wichtiger Schritt ist die
Konstruktion des Gitters $\tilde{L}$ im
Algorithmus~\ref{algo:improved_algorithm}, Zeile~13.
Wir ben"otigen hierf"ur ein Verfahren, das aus
einem Erzeugendensystem von h"ochstens $n+1$ Vektoren eine Basis von $\tilde{L}$
effizient berechnet. Weiterhin sollten diese Basisvektoren in irgendeinem
Sinn reduziert, d.h. hei"st m"oglichst kurz sein. Ohne diese Reduktion kann
es bei der Iteration von Basisberechnungen zu einer
Koeffizientenexplosion kommen, welche die Laufzeit der sp"ateren Basisberechnungen
au"serordentlich verschlechtert.  Dieses
Verhalten ist beispielsweise bei der 
Berechnung der Hermite-Normalform problematisch (siehe das Beispiel
von Hafner und McCurley in \cite{hafner89:_asymp}). Wir ziehen daher
die Verwendung einer Variante des LLL-Algorithmus vor, der bereits
w"ahrend der Basisberechnung 
reduziert.  Es sind exakte Schranken f"ur die L"ange von LLL-reduzierten
 Basisvektoren in
Beziehung zu den sukzessiven Minima des Gitters bekannt. F"ur die
maximale Laufzeit des
LLL-Algorithmus ist bei ganzzahligen Gittern eine Absch"atzung bekannt.
\begin{lem}\label{lemma:LLL}
  Sei $L$ ein ganzzahliges Gitter auf dem $n$-dimensionalen
  euklidischen Raum $E,$ welches von $n+1$ Vektoren, alle k"urzer als
  $B\in \BR,$ erzeugt wird. Dann kann eine (LLL-reduzierte) Basis
  $\{b_1,\cdots,b_n\}$ von $L$ in $O(n^4\log B)$ arithmetischen
  Operationen berechnet werden und es gilt f"ur $1\leqslant j\leqslant
  n,$ da"s
  \begin{equation}
    ||b_j||\leqslant 2^{(n-1)/2}\cdot \lambda_j(L)\leqslant
    2^{(n-1)/2}\cdot
    \lambda_n(L). \label{eq:lll_upper_bound_for_coeffs}
  \end{equation}
\end{lem}
\begin{proof}
  Nach \cite{buchmann89:_comput}, Theorem~3.1.\ l"a"st sich eine
  LLL-reduzierte Basis in $O(2n+2)^4\log B \subseteq O(n^4\log B)$
  arithmetischen Operationen berechnen. Die L"ange der Basisvektoren
  ist nach \cite{lenstra82:_factor}, Proposition~1.12  wie
  behauptet beschr"ankt.
\end{proof}
Nach diesen Vorbereitungen k"onnen wir die Laufzeit von
Algorithmus~\ref{algo:improved_algorithm} nach oben absch"atzen.
\begin{satz}\label{theorem:runtime_of_algorith_1}
  Sei  $L$ ein ganzzahliges Gitter auf dem $n$-dimensionalen
  euklidischen Raum $E, B\in\BR,$ so da"s
  \begin{math}
    S:=\{v_1,\cdots,v_s\} = \{v\in V\setminus 0 \mid ||v||
    \leqslant B\}
  \end{math}
  das Gitter  $L$ erzeugt. Dann k"onnen wir eine orthogonale Zerlegung
  \begin{math}
    L= L_1\oplus\ldots\oplus L_r
  \end{math}
  von $L$ in unzerlegbare Teilgitter $L_i,1\leqslant i\leqslant r$ in
  h"ochstens   $O(n^6\log^2 (nB) + sn^2)$ arithmetischen
  Operationen berechnen.
\end{satz}
\begin{proof}
  Wir analysieren die Laufzeit von
  Algorithmus~\ref{algo:improved_algorithm}. Dazu benutzen wir die
  in Abschnitt~\ref{sec:improved_algo_data_struct} beschriebenen
  Datenstrukturen.

  Mittels Radix-Sort (\cite{knuth74}) sortieren wir alle Vektoren in
  $S$ ihrer Norm nach und f"ugen sie in eine Liste ein. Die Norm
  berechnen wir in $O(n)$ Operationen, so da"s der gesamte
  Sortiervorgang in $O(sn)$ Operationen erfolgen kann. 
  
  Seien nun in Zeile 7 bei Eintritt in die Schleife $v_{11}, \cdots,v_{kr_k}, w_1,\cdots,w_{n-l}$ und $A$
  bereits per Induktion berechnet.  In Zeile~8 wird ein  $v\in S$ mit
  minimaler Norm entnommen  und gepr"uft, ob $v\in L_1\oplus\ldots\oplus L_k$,
  d.h., ob
die ersten $n-l$ Koeffizienten von
  $Av$ in $\BZ$ liegen und die restlichen $l$ Koeffizienten gleich $0$
  sind. Dies kostet f"ur alle $v\in S$ zusammengenommen $O(sn^2)$
  Operationen. Falls diese Bedingung erf"ullt ist, dann liegt $v$ in
  dem bereits erzeugten Untergitter und es wird kein Update 
der  Datenstrukturen ben"otigt. Nach
Satz~\ref{proposition:bounded_chain}~ii) werden  f"ur
  ganz $S$ in keinem Fall mehr als
  $\frac{1}{2}\log(n!(\frac{B}{M})^n)+n$ Updatevorg"ange benötigt. Bei
  jedem 
  Update wird eine LLL-reduzierte Gitterbasis aus einem
  Erzeugendensystem mit h"ochstens $n+1$ Vektoren berechnet. Diese
  Vektoren sind entweder aus $S$ oder Output eines fr"uheren LLLs, also
  LLL-reduziert. Damit ist ihre Norm durch $2^{(n-1)/2}\lambda_n(L)\leqslant
  2^{(n-1)/2}B$ nach oben beschr"ankt. Nach Lemma~\ref{lemma:LLL}
  k"onnen wir eine Basis des erzeugten Gitters in h"ochstens
  $O(n^4\log (2^{(n-1)/2}B))= O(n^5+n^4\log B)$ Operationen
  berechnen. Die Matrix $A$ k"onnen wir beispielsweise mit dem
  Gau"sverfahren (siehe \cite{golub96:_matrix}) in
  $O(n^3)$ Operationen updaten. 
  
  In der Summe erhalten wir damit 
f"ur die Gesamtlaufzeit eine obere Schranke von
  $O(n^6\log(nB) + n^5 \log(nB)\log B + sn^2)$ arithmetischen
  Operationen. Daraus folgt die leicht schlechtere, aber in der
  Schreibweise übersichtlichere
 Schranke von $O(n^6\log^2 (nB) +sn^2)$
  arithmetischen Operationen.
\end{proof}
\begin{remark}\label{remark:minimalvektoren_ohne_basis}
  Conway und Sloane konstruierten 1995 ein Gitter, das aus seinen
  Minimalvektoren erzeugt wird, aber keine Basis aus
  Minimalvektoren besitzt \cite{conway95}. Damit ist gezeigt, da"s
  es nicht m"oglich ist, das Anwachsen der Normen der Basisvektoren in
  der neuen
  Basis gegen"uber dem ursprünglichen Erzeugendensystem zu
  verhindern. Dennoch ist 
  die Absch"atzung "uber die L"ange von LLL-reduzierten Vektoren in Gleichung (\ref{eq:lll_upper_bound_for_coeffs}) "au"serst
  pessimistisch. In der Praxis findet der LLL-Algorithmus deutlich
  k"urzere Basisvektoren, insbesondere wenn das Erzeugendensystem
  bereits aus relativ kurzen Vektoren besteht. 
\end{remark}
\section{Gitterbasen aus großen
  Erzeugendensystemen}
\label{sec:mmlll}

Ein in der Praxis wichtigeres und viel h"aufigeres Problem als die
Zerlegung von Gittern ist die Berechnung einer Gitterbasis aus einem
gro"sen Erzeugendensystem. Hierzu stellen wir die bekanntesten
Algorithmen samt einiger Modifikationen zusammen. Im folgenden sei $L$
ein ganzzahliges Gitter mit Minimum $M$ im $n$-dimensionalen
euklidischen Raum $E,$ das von $s$ Vektoren erzeugt wird, die alle nicht
l"anger als $B\in \BR$ sind. Weiter nehmen wir an, da"s $s>n$ gilt:

\subsection{Der BP-LLL}\label{sec:BP-LLL}
Die schon in Lemma~\ref{lemma:LLL} angesprochene Modifikation des
LLL-Algorithmus von Buchmann und Pohst (\cite{buchmann89:_comput}), im
folgenden BP-LLL genannt, wird in der Praxis selten eingesetzt. F"ur
theoretische Überlegungen ist dieses Verfahren wertvoll, weil dafür eine
obere Absch"atzung der Laufzeit von 
$O((s+n)^4\log B)\subseteq O(s^4\log B)$ arithmetischen
Operationen existiert.

\subsection{Der MLLL}
\label{sec:MLLL}
Die von Pohst in \cite{pohst87:_lll} beschriebene Modifikation des LLL ist
unseres Wissens der
meistverwendete Algorithmus zur Basisbestimmung.  Leider ist die theoretische
Laufzeit des MLLL gegenw"artig nicht bekannt.  In der Praxis ist er jedoch in seinem
Laufzeitverhalten f"ur kleine Erzeugendensysteme (d.\ h.\ $s\approx n$) dem LLL
"ahnlich.% 
\subsection{Der serielle BP-LLL}
\label{sec:serieller_BP-LLL}
Anstatt den BP-LLL auf das gesamte Erzeugendensystem mit $s$ Vektoren
anzuwenden, benutzen wir ihn sukzessiv f"ur Erzeugendensysteme bestehend
aus $n+1$ Vektoren. Die L"ange der Vektoren in den dabei konstruierten Basen
lassen sich mit Lemma~\ref{lemma:LLL} durch
$2^{(n-1)/2}\lambda_n(L)\leqslant 2^{(n-1)/2}B$ nach oben
absch"atzen. Damit ist die Laufzeit des seriellen BP-LLL in
$O(s(2n+2)^4\log( 2^{(n-1)/2}B))\subseteq O(sn^4(n+\log B)).$

\subsection{Der serielle BP-LLL mit Aussieben}
\label{sec:LLL-Aussieben}

Wir k"onnen den Algorithmus~\ref{algo:improved_algorithm}
vereinfachen, um lediglich eine Gitterbasis statt einer Zerlegung des
Gitters  zu berechnen. Der daraus entststehende
Algorithmus~\ref{algo:mmlll} (s.~S.~\pageref{algo:mmlll}) benutzt eine
Unterfunktion \texttt{construct\_basis}, die eine Basis aus h"ochstens
$n+1$ Erzeugern berechnet. Um die Laufzeit mit
Satz~\ref{theorem:mmlll} abzusch"atzen, benutzen wir hierf"ur den
BP-LLL; man beachte aber auch hier die Bemerkung in~\ref{sec:MLLL}.

\begin{algorithm}[h]
\caption{Berechnung der Basis von $L$ aus einem Erzeugendensystem $S$}
\label{algo:mmlll}
\begin{algorithmic}[1]
\STATE \emph{Input:} Ein Erzeugendensystem $S=\{v_1,\cdots,v_s\}$ von
$L$ mit $||v_i||\leqslant B.$ 
\STATE \emph{Output:} Eine Basis $\{b_1,\cdots,b_n\}$ von $L.$

  \STATE {\sl 1.\ //Initialisierung.}
  \STATE W"ahle $v \in S, S \leftarrow S\backslash\{v\}, b_1\leftarrow v$

  \STATE {\sl 2.\ //F"uge Vektoren $v$ aus $S$ sukzessive hinzu.}
    \WHILE {$S \neq \emptyset$} 
  \STATE W"ahle $v \in S$ und setze $S \leftarrow   S\backslash\{v\}$.  
 \IF {$v \not\in \Zspan{b_1,\cdots,b_k}$}
 \STATE $\{b_1,\cdots,b_k\} \leftarrow \text{  \texttt{construct\_basis}}
 (b_1,\cdots,b_k,v)$
 \ENDIF \ENDWHILE
\end{algorithmic}
\end{algorithm}

\begin{satz}\label{theorem:mmlll}
  Sei  $L$ ein ganzzahliges Gitter auf dem $n$-dimensionalen
  euklidischen Raum $E,$ das von $s$ Vektoren, alle nicht l"anger als
  $B\in\BR,$ erzeugt wird. Dann k"onnen wir eine Basis von $L$ in
  h"ochstens   $O(n^6\log^2 (nB) + sn^2)$ arithmetischen
  Operationen berechnen.
\end{satz}
\begin{proof}
  Wir benutzen f"ur die Formulierung von Algorithmus~\ref{algo:mmlll} dieselben
  Datenstrukturen wie sie in Abschnitt~\ref{sec:improved_algo_data_struct}
  beschrieben sind und verwenden den BP-LLL f"ur die Unterfunktion
  \texttt{construct\_basis}. Dann berechnet der folgende
  Algorithmus~\ref{algo:mmlll} 
  eine (LLL-reduzierte) Basis von $L.$ Die Absch"atzung der
  Laufzeit erfolgt vollst"andig analog zum Beweis von
  Satz~\ref{theorem:runtime_of_algorith_1}. 
\end{proof}
\begin{remark}  
  Im Vergleich zum seriellen BP-LLL werden in
  Algorithmus~\ref{algo:mmlll} f"ur gro"se $s$ die meisten   LLLs
  durch eine einzige Matrix-Vektor-Multiplikation ersetzt. Daf"ur
  werden die 
  im Vergleich zu $s$ relativ wenige Updates mit der Berechnung der
  Inversen einer $n\times n$-Matrix verteuert.
\end{remark}
\section{Beispiele und Beobachtungen}
\label{sec:experimental_results}

\begin{bsp}\label{bsp:all_vectors_needed}
  Sei $E=\BR^1$ und $\{p_1,\cdots,p_k\}$ eine Menge von paarweise
  verschiedenen Primzahlen, deren Produkt gleich $t$ ist. Dann folgt
  mit dem chinesischen Restesatz $\langle
  \frac{1}{p_1},\cdots,\frac{1}{p_k} \rangle_\BZ = \frac{1}{t}\BZ.$ In
  dem Erzeugendensystem liegt kein Vektor im Erzeugnis der anderen.
  Das Erzeugendensystem ist also \glqq groß\grqq{} und minimal, aber
  das daraus erzeugte 
  Gitter besitzt das im Vergleich zu den $p_i$
  relativ kleine Minimum~$\frac{1}{t}.$
\end{bsp}
Bei diesem Beispiel handelt es sich um einen atypischen Extremfall. Viele
Gitter, die in der Theorie der quadratischen Formen behandelt werden, etwa
ganzzahlige Gitter, Gitter mit hoher Packungsdichte oder durch algebraische
Codes konstruierte Gitter, besitzen eher gro"se Minima. In disem Fall wird
der "Uberdeckungssatz~\ref{satz:covering_theorem} praktisch relevant und
besagt, da"s minimale Erzeugendensysteme solcher Gitter relativ klein sind.

Im folgenden untersuchen wir die praktische Anwendbarkeit von
Satz~\ref{satz:covering_theorem} und Algorithmus~\ref{algo:mmlll}. Als Kern
von Algorithmus~\ref{algo:mmlll} benutzen wir in unseren Versuchen den MLLL
von Pohst f"ur ganzzahlige Gitter. Bei nicht ganzzahligen Gittern
verkompliziert sich die Situation lediglich bei der Vorbereitung des
Inputs f"ur unseren 
Algorithmus. Ist das Gitter zwar nicht ganzzahlig, aber rational, so
k"onnen wir es durch geeignetes Skalieren auf ein ganzzahliges Gitter
transformieren. Handelt es sich um ein echt reelles  Gitter, so werden wir
zun"achst eine Approximation auf ein rationales Gitter vornehmen. Dieses
Verfahren wird beispielsweise in \cite{buchmann89:_comput} benutzt und seine
Auswirkungen auf den LLL analysiert. Um die Problematik der mit einer
Approximation verbundenen Rundungsfehler hier zu vermeiden, betrachten
wir nur ganzzahlige Gitter.

\begin{figure}[htbp]
  \includegraphics[width=0.6\textwidth,angle=270]{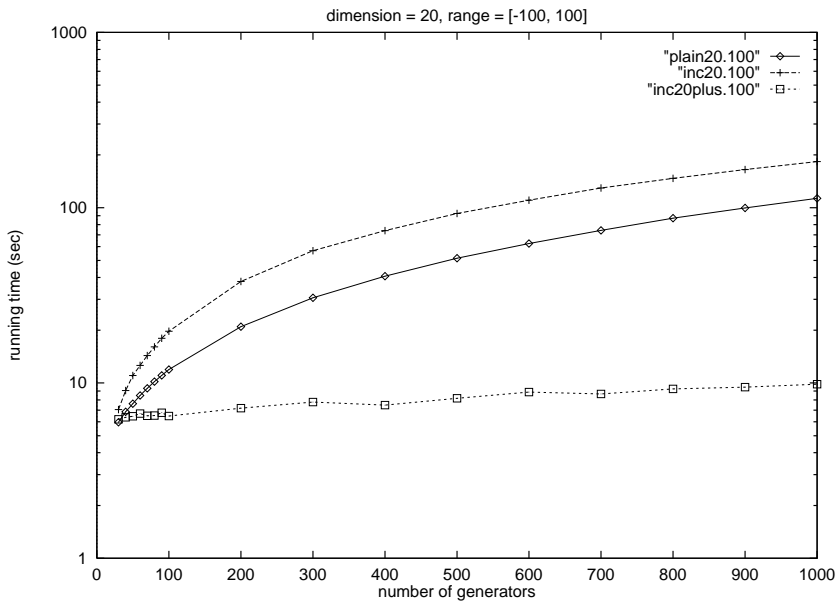}\vspace{4mm}

{\small $\lozenge$ MLLL, \hspace{1em} $+$ incremental
  MLLL,\hspace{1em}
  $\square$ Algorithm~\ref{algo:mmlll}} 
  \caption{Laufzeiten in  Dimension $20$ mit $K=100$}
  \label{fig:dim20_100}
\end{figure}

\begin{figure}[htbp]
\includegraphics[width=0.6\textwidth,angle=270]{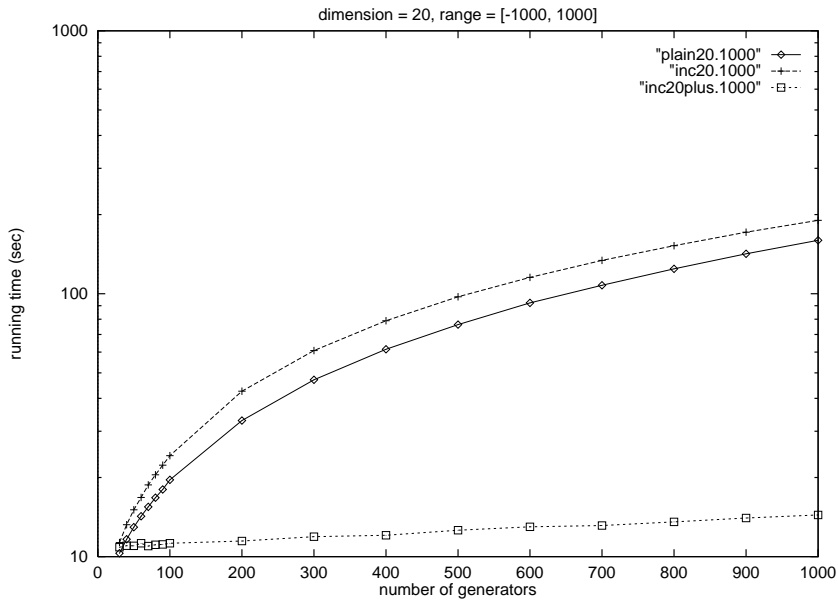}\vspace{4mm}

{\small $\lozenge$ MLLL, \hspace{1em} $+$ incremental
  MLLL,\hspace{1em}
  $\square$ Algorithm~\ref{algo:mmlll}} 
  \caption{Laufzeiten in  Dimension $20$ mit $K=1000$}
  \label{fig:dim20_1000}
\end{figure}

Lenstra, Lenstra und Lov\'asz beweisen in \cite{lenstra82:_factor}  f"ur die Laufzeit des LLLs
eine asymptotische obere Schranke von $O(n^4\log B)$ arithmetischen
Operationen. Experimente zeigen, da"s diese
Schranke pessimistisch ist; dies gilt ebenfalls f"ur die oben
erw"ahnten Modifikationen des LLL wie den MLLL. Bei sehr vielen
Inputvektoren wird jedoch auch der MLLL merkbar langsam. Die naive Idee,
sukzessive den MLLL auf $n+1$ Vektoren anzuwenden, führt zu keiner 
Verbesserung.  Wir nennen diese Variante den \emph{inkrementellen
  MLLL}. Die Abbildungen \ref{fig:dim20_100} und \ref{fig:dim20_1000}
zeigen die Laufzeiten vom MLLL, dem inkrementellen MLLL und
Algorithmus~\ref{algo:mmlll} im Vergleich. Wir wenden die Algorithmen
auf Erzeugendensysteme mit Vektoren an, deren Koeffizienten
ganzzahlig, pseudo-zuf"allig und vom Betrag kleiner als eine Konstante
$K$ sind. Insbesondere ist die Norm der Vektoren durch $B= \sqrt{n}K$
beschr"ankt. 

\bigskip

Die Auswertung der Messungen in den Abbildungen \ref{fig:dim20_100}
und  \ref{fig:dim20_1000} zeigt, 
da"s die Anzahl der Updateoperationen in 
Algorithmus~\ref{algo:mmlll} sehr klein in Beziehung zur Gr"o"se des
Erzeugendensystems ist. Typische Determinanten der Anfangsgitter von
 vollem Rang sind in der Gr"o"senordnung von $2^{100},$ die
Determinanten der erzeugten Gitter sind sehr klein. Die Determinanten
der Anfangsgitter enthalten oft sehr gro"se Primteiler, so da"s die
Anzahl der Updateoperationen sehr viel kleiner ist als die
in Korollar~\ref{kor:number_of_update_steps} angegebene Schranke, bei der die
vorkommenden Primteiler nach unten auf $2$ abgesch"atzt werden.

Die Abbildungen \ref{fig:dim20_100} und \ref{fig:dim20_1000} stammen
von Frank Vallentin, der auch das Testprogramm entwickelt hat.  Dabei
wurde als Kern Victor Shoups NTL\footnote{URL:
  http://www.shoup.net/ntl/} Bibliothek benutzt, aus
welcher die rationale Arithmetik und der MLLL stammen.

\section{Implementierungshinweise}
\label{sec:mlll_performance}

Die oben beschriebenen Verfahren führen, wie an den Beispielen gezeigt,
zu in der Praxis effizienten Implementierungen mit beweisbaren
Komplexit"atsabsch"atzungen. Dennoch sind mit einfachen Heuristiken und
Methoden in vielen F"allen weitere erhebliche Laufzeitgewinne zu erzielen. Die
Liste der folgenden Hinweise sollte nur als Wegweiser zu einer
Implementation von schnellen und robusten Bib\-li\-otheks\-funk\-ti\-on\-en
dienen.

\subsection{Das Zerlegungsproblem}
\label{sec:mlll_perf_zerl}

Beim Zerlegungsproblem ist eine Basis oder Grammatrix des Gitters
bereits bekannt. Daraus kann man unmittelbar die Invarianten Rang und
Determinante des 
Gitters sehr leicht berechnen. Diese können als Abbruchkriterium
im Schleifendurchlauf von Algorithmus~\ref{algo:improved_algorithm}
verwendet werden.
\begin{enumerate}
\item Wenn im Anschluß an Zeile 15 in
  Algorithmus~\ref{algo:improved_algorithm} das Gitter $L_1$
  denselben Rang wie das Gitter $L$ besitzt, dann kann es nur ein
  triviales orthogonales Komplement besitzen, also ist bereits auch das
  Gitter $L$ unzerlegbar und der Algorithmus kann beendet werden.
\item  Wenn im Anschluß an Zeile 15 in
  Algorithmus~\ref{algo:improved_algorithm} die Determinante
  des Gitters $L_1\oplus \ldots\oplus L_K$ gleich der Determinate von
  $L$ ist, dann 
  gilt $L_1\oplus \ldots\oplus L_K =L,$ und wir haben eine Zerlegung von $L$
  gefunden, wenn auch noch nicht alle Vektoren in $S$ abgearbeitet sind.
\end{enumerate}

Ist ein Gitter zerlegbar, so kann diese Eigenschaft gelegentlich schon an einer
LLL-reduzierten Grammatrix des Gitters abgelesen werden. Wenn
beispielsweise alle sukzessiven Minima nahe beieinander liegen, dann
sind zerlegbare Vektoren relativ lang im Vergleich zu unzerlegbaren
Vektoren. Der LLL-Algorithmus wird dann mit großer Wahrscheinlichkeit
bereits eine Basis aus unzerlegbaren Vektoren finden. Die zu dieser
Basis gehörige Grammatrix besitzt dann (bis auf eine
Permutation der Basisvektoren) die Gestalt einer Blockmatrix.

\subsection{Die Basisberechnung aus einem großen
  Erzeugendensystem}\label{sec:mmlll_perf_basis} 
Der Algorithmus~\ref{algo:mmlll} ist leicht parallelisierbar. Eine
Erzeugendenmenge kann in kleinere, etwa gleich große Teilmengen
zerlegt werden. Die Basen der aus den Teilmengen erzeugten Gitter
können jeweils mit dem Algorithmus~\ref{algo:mmlll} nach dem "`divide and
conquer"'-Verfahren parallel berechnet
werden. In einem letzten Schritt werden die Basisvektoren der
Teilgitter wieder zu einem Erzeugendensystem $S$ zusammengefügt,
welches wiederum als Input für Algorithmus~\ref{algo:mmlll} dient. Diese
Vorgehensweise empfiehlt sich allerdings nur für sehr große
Erzeugendensysteme. In der Regel ist selbst eine h"aufige Ausführung der
Lookup-Funktion in Zeile~8 schneller als der
zus"atzliche Aufwand durch wenige neue MLLLs.

\chapter{Der Quantizer  der Vorono"izelle eines Gitters}
\label{sec:sec:quantizer}
Quadratische Formen und ihre Gitter finden in einem Gebiet der
Informationstheorie wichtige Anwendungen. Bei der "Ubertragung von
Informationen "uber ein st"orungsanf"alliges Medium ist es von eminenter
Bedeutung, auf der Empf"angerseite die gesendeten Informationen korrekt und
vollst"andig wieder herstellen zu k"onnen. Die grunds"atzliche Idee zur
Fehlerkorrektur besteht in der Erh"ohung der Redundanz der Informationen, um
kleinere Fehl"ubertragungen korrigieren zu k"onnen. Eine solche M"oglichkeit
ist die Verwendung fehlerkorrigierender Codes. Viele dieser Codes mit guten
Korrektureigenschaften k"onnen aus Gittern konstruiert werden. Diese Methoden
greifen aber erst, wenn die Informationen bereits als codierte, digitale Daten
vorliegen. Bei der "Ubertragung durch Kabel oder Funkwellen erreichen die
Informationen den Empf"anger zun"achst analog und in der Regel gestört. Ein
Analog-Digital-Wandler setzt die analogen Signale in digitale um. Schon in
diesem Schritt werden kleinere St"orungen korrigiert. Hierf"ur werden die
gesendeten Daten nur aus einem kleinen, fest definierten Zeichenvorrat, dem
sogenannten Code Book, gebildet. Ein empfangenes, gest"ortes Datum wird
mittels eines Quantizers, manchmal auch Approximation oder Compressor genannt,
auf das "`n"achstliegende"'  Zeichen im Code Book abgebildet. Bereits diese
informale Beschreibung l"a"st die Beziehung zwischen Quantizern und
Vorono"izellen erahnen. Eine Methode, Quantizer mit guten
Korrektureigenschaften zu konstruieren, basiert auf Gittern im euklidischen Raum.

\begin{definition}
  Sei ${P}$ eine diskrete Punktmenge in $E.$ Eine Funktion
  $q:E \rightarrow {P}$ heißt Quantizer, wenn für jedes
  $v\in E$ gilt, daß 
  \begin{math}
    ||v-q(v)|| \leqslant ||x-p|| \text{ für alle } p\in{P}.
  \end{math}
\end{definition}
Ein Quantizer ordnet  also jedem Punkt im euklidischen Raum den
n"achstliegenden Punkt oder, falls dieser nicht eindeutig ist, einen der
n"achstliegenden Punkte in ${P}$ zu. 

Ein Quantizer ist also im wesentlichen bereits durch die
Angabe seiner Wertemenge charakterisiert. Lediglich für Punkte in $E,$
die sich genau gleich weit von zwei Punkten in $P$ befinden, ist der
Funktionswert f"ur zwei Quantizer von $P$ nicht notwendigerweise gleich. 

Man kann sich die Funktion eines Quantizers als hochdimensionalen
Rundungsprozeß oder auch als eine verlustbehaftete Datenkompression
vorstellen. Für die Qualit"at eines Quantizers k"onnen je nach Anwendung
verschiedene Ma"se dienen. Ein universelles und weit verbreitetes ist der
durchschnittliche quadratische Fehler je Dimension
  \begin{equation}
  G_q =\frac{1}{n}\int_{E} ||v-q(v)||^2\cdot p(v)\,dv,\label{eq:quantizer}
\end{equation}
wobei die Funktion $p$ die Wahrscheinschlichkeitsdichte der Datenquelle
beschreibt (\cite{bible}, S.~58).

Wenn man für $P$ die Menge der Punkte eines Gitters $L$ w"ahlt, so l"aßt
sich die Konstruktion von gitterbasierten Quantizern elegant mit dem
geometrischen Konzept der Vorono"izelle eines Gitters beschreiben. Im
folgenden Abschnitt werden einige wichtige Definitionen und
Eigenschaften der Vorono"izellen von Gittern zusammengefaßt.

Alle Gitter $L$ in diesem Kapitel haben vollen Rang $n,$ d.h., sie
spannen $E$ als $\BR$-Vektorraum auf.

\section{Die Vorono"izelle eines Gitters}
\label{sec:voronoizellen}

Zun"achst wird der zentrale Begriff des Vorono"ibereichs eines Punktes
eingeführt. 

\begin{definition}\label{def:voronoibereich_voronoizelle}
  Sei $P$ eine diskrete Teilmenge von $E.$ F"ur ein $v\in P$ hei"st
  die Menge $$V_v(P) = \{ x\in E \mid ||x-v|| \leqslant ||x-w|| \text{
    f"ur alle } w\in P\}$$
  der \emph{Vorono"ibereich von $v$ in $P.$}
  
  Ist $v\in L$ ein Gitterpunkt, so nennt man $V_v(L)$ die \emph{Vorono"izelle
    von $v$ in $L.$} Wenn keine Verwechselung m"oglich ist, dann schreibt man
  auch kurz $V_v.$ 
\end{definition}
Offensichtlich gilt $q^{-1}(v)\subseteq V_v(P)$ für ein $v\in P.$ Im
allgemeinen gilt jedoch keine Gleichheit, weil es Punkte $x\in E$
gibt, die sich in mehr als einer Vorono"izelle befinden.

Die Punkte in $E,$ die sich am weitesten von Gitterpunkten entfernt
befinden, kennzeichnen eine wichtige Invariante des Gitters, den
Überdeckungsradius. 
\begin{definition}\label{def:covering_radius}
  F"ur eine Teilmenge $P\subseteq E$ hei"st die Zahl
\begin{displaymath}
  R(P) = \sup_{x\in E} \inf_{v\in P} ||x-v|| \in \mathbb{R}^{\geqslant
  0}\cup\{\infty\}
\end{displaymath}
der \emph{"Uberdeckungsradius} von $P.$
\end{definition}

 R(L) ist der gr"o"stm"ogliche Abstand,
den ein Punkt in $E$ von einem Gitterpunkt haben kann. Das Auf\/finden von
Punkten mit dieser Eigenschaft, den sogenannten \emph{tiefen L"ochern}\label{deep_hole}  ist hilfreich zur Konstruktion von Gittern mit
hoher Packungsdichte (s.~\cite{bible} Chapt.~6).

 In diesem Kapitel werden einige geometrischen Eigenschaften von
Vorono"izellen vorgestellt. 

\section{Die Vorono"izelle als Polytop}
\label{sec:Voronoizelle_ist_ein_Polytop}
Das  erste Ziel ist es, die Vorono"izelle $V_0$ des Ursprungs $0$ eines
Gitters  als Polytop zu beschreiben und ihre Facetten zu bestimmen.

In diesem Kapitel tritt der Aspekt der definierenden quadratischen Form in den
Hintergrund. Zur Vereinfachung wird hier lediglich das euklidische
Skalarprodukt $(v,v)$ benutzt, das wegen der zahlreichen Klammern in einigen
Formeln auch als $v\cdot v$ geschrieben wird.  Falls $v\in E\setminus 0,$
schreiben wir $H(v)=\{x\in E \,|\, x\cdot v = v\cdot v\}$ f"ur die Hyperebene
mit dem Normalenvektor $v$ und $H_v = \{x\in E \,|\, x\cdot v \leqslant v\cdot
v\}$ f"ur den dadurch bestimmten abgeschlossenen Halbraum, der die $0$
enth"alt.

Wenn $P$ eine diskrete Menge mit $0\in P$ ist, dann kann man in dieser
Notation den Vorono"ibereich von $0$ als Durchschnitt von Halbr"aumen
beschreiben.

\begin{satz}
  Sei $P$ eine diskrete Menge und $0\in P.$ Dann ist
  \begin{displaymath}\label{v_0_ist_durchschnitt}
    V_0 = \bigcap_{v\in P\setminus \{0\}} H_{\frac{1}{2} v}.
  \end{displaymath}
\end{satz}
\begin{proof} Sei $0\neq v\in P$ und $x\in E.$ Dann gilt
  \begin{displaymath}
    \begin{array}{rl}
    &||x-0|| \leqslant ||x-v|| \\
    \Longleftrightarrow &x\cdot x\leqslant
    (x-v)\cdot (x-v) = x\cdot x -2 x\cdot v + v\cdot v\\
    \Longleftrightarrow& x\cdot v \leqslant \frac{1}{2} v\cdot v\\
    \Longleftrightarrow& x\cdot \frac{1}{2}v \leqslant
    \frac{1}{2} v\cdot \frac{1}{2}v\\
    \Longleftrightarrow& x\in H_{\frac{1}{2}v}.\label{}
  \end{array}
\end{displaymath}
 Damit folgt
\begin{displaymath}
  x\in V_0 \quad \Longleftrightarrow\quad \forall_{0\neq v\in P}
  ||x-0|| \leqslant ||x-v|| \quad \Longleftrightarrow\quad x\in
  \bigcap_{0\neq v\in P} H_{\frac{1}{2} v}.
\end{displaymath}
\end{proof}

Im Fall, daß $P$ ein Gitter ist, folgt 
aus der Definition \ref{def:voronoibereich_voronoizelle} und der
Eigenschaft $P=P-v$ f"ur alle $v\in P$ unmittelbar das
\begin{lem}\label{lem:voronoizelle_ist_invariant}
  Sei $L$ ein Gitter auf $E.$
  \begin{enumerate}
  \item F"ur $v\in L$ gilt $V_v = v+V_0.$
  \item F"ur $\sigma\in O(L)$ gilt $V_0\sigma = V_0$, inbesondere ist
    $-V_0= V_0.$
  \end{enumerate}
\end{lem}

Die Vorono"izelle eines beliebigen Gitterpunktes ist ein Translat von $V_0.$
Deshalb schreiben wir auch "`die"' Vorono"izelle von $L,$ wenn wir nicht eine
ausgezeichnete Zelle meinen.

\begin{lem}\label{lem:voronoipflasterung_ist_face2face}
  Es seien $v,w\in V_w,\, x\in V_w$ und $v\in L\setminus \{w\}$ mit $||x-w|| =
  ||x-v||,$ 
  dann ist auch $x\in V_v.$
\end{lem}
\begin{proof}
  Sei $z\in L\setminus \{v, w\}$. Aus $x\in V_w$ folgt,
  da"s $||z-x||\geqslant||x-w||=||x-v||.$ Also ist $x\in V_v.$
\end{proof}

\begin{remark}\label{satz:face2face}
  Aus diesem Lemma folgt u.a. der wichtige Satz, da"s $\CT = \{ V_v\mid v\in
  L\}$ eine Seite-an-Seite Pflasterung des Raumes $E$ ist.
\end{remark}

\begin{lem}
  Sei $C\subseteq E$ konvex und abgeschlossen.  Die Menge $C$ ist
  genau dann beschr"ankt, wenn sie keinen Strahl enth"alt.
\end{lem}
\begin{proof}
  Wenn $C$ einen Strahl enth"alt, dann ist $C$ offenbar nicht
  beschr"ankt. Sei $C$ nun nicht beschr"ankt und ohne Beschr"ankung der
  Allgemeinheit  sei die Sph"are mit Radius $1$ um den Nullpunkt $S_0(1)= \{x\in E\, |\, ||x||=1\} \subseteq C.$ Da die Menge $C$ nicht beschr"ankt ist,
  enth"alt sie eine Folge $(x_n)_n \subseteq C$ mit $||x_n||>n.$ Wegen der
  Kompaktheit von $S_0(1)$  enth"alt die Folge
  $(\frac{x_n}{||x_n||})_n \subseteq S_0(1)$ nach dem
  Satz von Bolzano-Weierstra"s
eine konvergente
  Teilfolge $(\frac{x_{n_k}}{||x_{n_k}||})_{n_k},$ deren Limes wir mit
  $x$ bezeichnen. Wir zeigen nun, da"s $C$ den Strahl $ R^+ x$
  enth"alt. Sei $\alpha\in R^+.$ F"ur fast alle $k$ ist
  $x_{n_k}>\alpha$ und wegen der Konvexit"at von $C$ dann auch
  $\frac{\alpha}{||x_{n_k}||}x_{n_k} \in \left[ 0,x_{n_k}\right]
  \subseteq C\cap S_0(\alpha) .$ Aus der Abgeschlossenheit von $C$ folgt, 
da"s $\lim_k
  \frac{\alpha}{||x_{n_k}||}x_{n_k}= \alpha x \in C$ und damit die
  Behauptung.
\end{proof}
\begin{lem}\label{lem:v_0_ist_beschraenkt}
  $V_0(L)$ ist beschr"ankt.
\end{lem}
\begin{proof}
  Seien $w,x\in E$ mit $w+\alpha x\in V_0$ f"ur alle $\alpha\in
  \BR^{\geqslant 0}$ und $S_x=\{v\in L\mid x\cdot v \geqslant 0\}.$
  F"ur alle $\alpha\in\BR^{\geqslant 0}$ und $v\in S_X$ gilt, da"s
  $w+\alpha x\in V_0 \Leftrightarrow (w+\alpha x)\cdot \frac{1}{2}v
  \leqslant \frac{1}{2}v\cdot \frac{1}{2}v.$ Damit folgt $x\cdot v =
  0.  $ $L$ ist ein Gitter auf $E$ und somit ist $\dim \langle
  S_x\rangle_\BR = \dim E$ und folglich $x=0.$ $V_0$ ist konvex und
  enth"alt keinen Strahl, womit ist gezeigt, da"s $V_0$ beschr"ankt
  ist.
\end{proof}

\section{Vorono"irelevante Vektoren}
\label{sec:voronoirelevante_vektoren}

\begin{definition}\label{def:voronoirelevant}
Sei $v\in L\setminus 0.$
\begin{enumerate}
\item  Wir nennen $v$ einen \emph{H"ullenvektor,} wenn
  $H(\frac{1}{2}v)\cap V_0\neq\emptyset,$ also seine Hyperebene die
  Vorono"izelle ber"uhrt.
\item $v$ hei"st \emph{vorono"irelevant} oder kurz \emph{relevant,}
  wenn $H(\frac{1}{2}v)$ einen $n-1$-di\-men\-sio\-na\-len (affinen) Schnitt
  mit $V_0$ besitzt, d.~h., $H(\frac{1}{2}v)$ eine Wand von $V_0$
  enth"alt \cite{bible}. Ist $v$ relevant, so schreiben wir
  \emph{$F_v$} f"ur die Wand $H(\frac{1}{2}v)\cap V_0.$
\end{enumerate}
\end{definition}

Die Menge der vorono"irelevanten Vektoren gen"ugt, um $V_0$ zu
konstruieren.

\begin{satz}\label{satz:waende_sind_zentralsymmetrisch}
  Die W"ande $F_v$ von $V_0$ sind zentralsymmetrisch mit Mittelpunkt
  $\frac{1}{2}v.$
\end{satz}
\begin{proof}
  Sei $t\in F_v.$ Wie in \ref{lem:voronoizelle_ist_invariant} gezeigt, liegt
  $-t+v$ in $V_v$ und wegen $(-t+v)\cdot \frac{1}{2}v= \frac{1}{2}v \cdot
  \frac{1}{2}v$ auch in $H(\frac{1}{2}v),$ also mit
  Lemma~\ref{lem:voronoipflasterung_ist_face2face} in $V_0\cap H(\frac{1}{2}v)
  = F_v.$ Damit ist auch die Konvexkombination $\frac{1}{2}t +
  \frac{1}{2}(-t+v) = \frac{1}{2}v$ in $F_v$ und $F_v$ ist zentralsymmetrisch
  bez"uglich $\frac{1}{2}v.$
\end{proof}
Wir wissen nun, da"s $V_0$ als Polytop ein endlicher Durchschnitt von
Halbr"aumen ist.  Bereits 1908 hat Vorono"i eine Kennzeichnung
dieser Halbr"aume angegeben:
\begin{satz}[\cite{voronoi08:_nouvel}]\label{voronoi_relevante_vektoren}
  Ein Vektor $0\neq v\in L$ ist genau dann relevant, wenn $\pm v$ die
  einzigen k"urzesten Vektoren in der Nebenklasse $v + L/2L$ sind. 
\end{satz}
\begin{proof}
  Wegen der Wichtigkeit dieses Satzes und der K"urze des Argumentes geben wir
  hier den Beweis aus \cite{bible}, S.~475\footnote{In der ersten Auf\/lage
    fehlt wegen eines Druckfehlers die Formulierung dieses Satzes. In der
    zweiten Auf\/lage ist der Satz enthalten, aber der Korrekturhinweis im
    Vorwort, S.~xxviii, verweist auf die falsche Seite.} an.
  
  Es seien $v,w\in L$ mit $v-w\in 2L, v\neq \pm w$ und $w\cdot w\leqslant
  v\cdot v.$ Dann sind $t=\frac{1}{2} (v+w)$ und $u=\frac{1}{2} (v-w)$ beide
  in $L\setminus 0$ und man verifiziert $x\cdot v=x\cdot(t+u)\leqslant
  \frac{1}{2}(t\cdot t + u\cdot u)\leqslant \frac{1}{4}(v\cdot v + w\cdot
  w)\leqslant \frac{1}{2}v\cdot v.$ Damit ist $v$ nicht relevant, weil
  $x\in H_t\cap H_u \Rightarrow x\in H_v.$
  
  Sei nun $v\in L$ nicht relevant und einer der k"urzesten Vektoren in der
  Nebenklasse $v+2L.$ Dann gibt es ein $w\in L\setminus 0$ mit $w\neq v$ und
  $\frac{1}{2}w\cdot \frac{1}{2}w\leqslant \frac{1}{2}v\cdot \frac{1}{2}w$.
  Wir erhalten
\begin{displaymath}
(v-2w)\cdot (v-2w)=v\cdot v - 4\,v\cdot w+4\,
w\cdot w \leqslant v\cdot v
\end{displaymath}
und haben mit $v-2w$ wegen $0\neq v-2w\neq \pm v$ einen weiteren k"urzesten
Vektor in $v+2L$ gefunden.
\end{proof}

An dieser Stelle sei nochmals auf das Korollar~\ref{kor:voronoirelevant} auf
Seite~\pageref{kor:voronoirelevant} verwiesen, das genau die spitz
unzerlegbaren Vektoren als vorono"irelevante identifiziert.

\begin{kor} Für ein Gitter $L$ gilt:
\begin{enumerate}
  \item Es gibt h"ochstens $2\cdot (2^n-1)$ relevante Vektoren. 
  \item $V_0(L)$ ist endlicher Durchschnitt von Halbr"aumen, also ein Polytop,
    insbesondere konvex.
  \item  $R(L)<\infty.$
\end{enumerate}
\end{kor}

\begin{remark}
  Bereits Vorono"i zeigte, da"s die Schranke von $2\cdot (2^n-1)$ relevanten
  Vektoren scharf ist \cite{voronoi08:_nouvel,voronoi09:_nouvel}.
\end{remark}

\section{Der Quantizer eines Gitters}
Im Falle einer Datenquelle mit uniform verteiltem Output und einem
Gitterquantizer
vereinfacht sich die 
Gleichung~(\ref{eq:quantizer}). Dies motiviert die folgende 
\begin{definition}[\cite{bible}, Chapt.~2.3, Chapt.~21.]\label{def:traegheitsmoment}
F"ur ein Gitter $L$ von vollem Rang, hei"st
\begin{equation}
  \label{eq:latticequantizer}
  G(L) := G_{V_0} = \frac{1}{n} (\det L)^{-\frac{n+2}{2n}}
  {\int_{V_0}x\cdot x\, 
  dx}\quad, 
\end{equation}
das \emph{normalisierte dimensionslose zweite  Tr"agheitsmoment} von $L.$
\end{definition}
Mathematisch exakt l"aßt sich diese Invariante von $L$ im allgemeinen nur
schwer errechnen. In \cite{bible}, Chapt.~21  findet sich die Herleitung
für einige wichtige Standardpolytope, unter ihnen auch das
Simplex. Letzteres wird dort benutzt, um mit Hilfe der Operation der
Weyl-Gruppe den Quantizer $G(L)$ für alle Wurzelgitter $L$ zu berechnen. 

Für Gitter mit weniger Strukturinformationen ist es praktikabel, das
Tr"agheitsmoment in Gleichung~(\ref{eq:latticequantizer}) durch
numerische Mon\-te-Car\-lo-In\-te\-gra\-tion zu approximieren. Es
ist jedoch nicht unmittelbar möglich, eine uniforme Verteilung von
Samples in $V_0$ zu erzeugen, weil die genaue Struktur von $V_0$ im
allgemeinen nicht bekannt ist. Man kann sich hier die
Translationseigenschaft von Vorono"izellen zunutze machen. Man
definiert eine Abbildung $V_0 \rightarrow [0,1)^n$ folgendermaßen. Für
eine feste Basis $B$ von $L$ bestimme man für jedes $y\in V_0$ den
Koordinatenvektor $u=B^{-1}$ bezüglich $B$ und bestimme
komponentenweise den nicht ganzzahligen Anteil $x=(x_1,\cdots,x_n)=
(u_1-\lfloor u_1\rfloor,\cdots, u_n-\lfloor u_n\rfloor).$ Mit Hilfe der
Translationseigenschaft von Vorono"izellen verifiziert man leicht, daß
diese Abbildung bijektiv ist. Somit kann man mit einer uniformen
Datenquelle in $[0,1)^n$ eine uniforme Verteilung in $V_0$
erzeugen.\footnote{Die zugrundeliegende geometrische Idee wird in
  \cite{conway82:_voron,conway82:_fast} beschrieben.} Die
Transformation mit $B$ muß im Integral in Gleichung~(\ref{eq:latticequantizer}) durch 
eine entsprechende   Substitution berücksichtigt werden.

Im folgenden Algorithmus~\ref{algo:latticequantizer}
 bezeichnet $\overline{S}$ das arithmetische Mittel
der  Abst"ande der Samples zu Gitterpunkten, $G(L)$ das
approximierte Tr"agheitsmoment, $s^2=\frac{1}{t-1}\sum_k
(S[k]-\overline{S})^2$  die  Streuung und
$\hat{\sigma}= \frac{s}{\sqrt{n}}$ den Sch"atzwert des Fehlers. Für den
korrekten Rückgabewert müssen $\hat{\sigma}$ und $G(L)$ 
jeweils mit $\frac{1}{n} (\det L)^{-\frac{1}{n}}$ skaliert werden.

\begin{algorithm}[h]
\caption{Numerische Approximation von $G(L)$}
\label{algo:latticequantizer}
\begin{algorithmic}[1]
\STATE \emph{Input:} Ein Gitter $L\subseteq E,$ gegeben durch eine
Basis $B$ und $t\in \BN,$ die Anzahl der zu berechnenden Samples. 

\STATE \emph{Output:} $G(L)$ und $\hat{\sigma},$ der Sch"atzwert des
Fehlers von $G(L).$

  \STATE {\sl // 1.\ Initialisierung.}
  \STATE Initialisiere ein Array $S[0...k]$ für die zu berechnenden  Samples.

  \STATE {\sl // 2.\ Berechnung der Samples.}
  \FOR{$k=0$ to $t$}
  \STATE W"ahle $x_k\in [0,1)^n$ randomisiert.
  \STATE  $y_k \leftarrow \closestvector(L,B\cdot x_k)$
  \STATE  $S[k] \leftarrow ||y_k-B\cdot x_k||^2$
  \ENDFOR
  \STATE {\sl // 3. Ergebnisaufbereitung.}
  \STATE ${\overline{S}} \leftarrow \frac{1}{n}  \sum_i^t S[i]$. {\sl
  // Arithmetisches Mittel über die $S[k]$}
  \STATE ${\sigma} \leftarrow   \sqrt{\frac{1}{t-1}
    \sum_i^t(\overline{S}-S[i])^2} $ 
  \STATE $\hat{\sigma} \leftarrow \frac{1}{n} (\det L)^{-\frac{1}{n}}\cdot
  \frac{\sigma}{\sqrt{t}}  $ 
  \STATE $G(L) \leftarrow   \frac{1}{n} (\det L)^{-\frac{1}{n}}\cdot
  \overline{S}$
  \STATE return $(G(L),\hat{\sigma})$
 \end{algorithmic}
\end{algorithm}

Algorithmus~\ref{algo:latticequantizer} kann leicht verallgemeinert
werden, um den mittleren quadratischen Fehler für Quantizer von
beliebigen Vereinigungen von Nebenklassen eines Gitters $L_C
\bigcup_{c\in C} L+c$ für $C=\{c_1,\cdots,c_m\} \subset E$ zu
berechnen. 
Hierzu bestimmt man nicht den Abstand von $B\cdot x_k$ zu $L,$
sondern jeweils den zu $L+c, c\in C$ und setzt $S[k]$ auf das Minimum
aller berechneten Abst"ande.  Für die Bestimmung von $G(L_C)$ und des
entsprechenden Sch"atzwertes des Fehler muß dem Umstand Rechnung getragen
werden, daß sich die Dichte der Gitterpunkte um den Faktor $m$
gegenüber der Dichte von $L$ vervielfacht. Dies wird durch
Korrekturfaktor $\frac{1}{|C|^2}$ bei der Determinante berücksichtigt.
Dem liegt die Identit"at $(\sum_{c\in C} \vol V_{c})^2 = \det L$
zugrunde, insbesondere ist es nicht notwendig, die exakten Volumina
$V_{c_1},\cdots, V_{c_m}$ zu kennen.  Somit erh"alt man Algorithmus
\ref{algo:latticecosetquantizer}.

\begin{algorithm}[h]
\caption{Numerische Approximation von $G(L_C)$}
\label{algo:latticecosetquantizer}
\begin{algorithmic}[1]
\STATE \emph{Input:} Ein Gitter $L\subseteq E,$ gegeben durch eine
Basis $B,$ eine Menge $C=\{c_1,\cdots,c_m\} \subset E$ und $t\in \BN,$
die Anzahl der zu berechnenden Samples.  

\STATE \emph{Output:} $G(L)$ und $\hat{\sigma},$ der Sch"atzwert des
Fehlers von $G(L)$

  \STATE {\sl // 1.\ Initialisierung.}
  \STATE Initialisiere ein Array $S[0...k]$ für die zu berechnenden  Samples.

  \STATE {\sl // 2.\ Berechnung der Samples.}
  \FOR {$k=0$ to $t$}
  \STATE W"ahle $x_k\in [0,1)^n$ randomisiert.
  \STATE  $y_k \leftarrow \closestvector(L,B\cdot x_k)$
  \STATE  $S[k] \leftarrow ||y_k-B\cdot x_k||^2$
  \FOR {$j=1$ to $m$}
     \STATE  $y_k \leftarrow \closestvector(L,B\cdot x_k-c_j)$
     \STATE  $S[k] \leftarrow \min(S[k], ||y_k-B\cdot x_k+c_j||^2)$
  \ENDFOR
  \ENDFOR
  \STATE {\sl // 3. Ergebnisaufbereitung.}
  \STATE ${\overline{S}} \leftarrow \frac{1}{n}  \sum_i^t S[i]$. {\sl
  // Arithmetisches Mittel über die $S[k]$}
  \STATE ${\sigma} \leftarrow   \sqrt{\frac{1}{t-1}
    \sum_i^t(\overline{S}-S[i])^2} $ 
  \STATE $\hat{\sigma} \leftarrow \frac{1}{n} (\det L/|C|^2)^{-\frac{1}{n}}\cdot
  \frac{\sigma}{\sqrt{t}}  $ 
  \STATE $G(L_c) \leftarrow   \frac{1}{n} (\det L/|C|^2)^{-\frac{1}{n}}\cdot
  \overline{S}$
  \STATE return $(G(L_c),\hat{\sigma})$
 \end{algorithmic}
\end{algorithm}

\section{Gitter mit kleinem Quantizer}
Das dimensionslose zweite normalisierte Tr"agheitsmoment eines
Quantizers $G(L)$ ist eine Invariante des Gitters $L,$ die bis vor
kurzem wenig untersucht wurde. Aus der Definition von $G(L)$ in
Gleichung~(\ref{eq:latticequantizer}) folgt unmittelbar, daß für ein
festes Volumen von $V_0(L)$ das Tr"agheitsmoment minimal ist, wenn
${\int_{V_0}x\cdot x\, dx}$ minimal für alle Gitter $L$
ist. Nach dem Isoperimeterprinzip bes"aßen $n$-dimensionale Kugeln das
kleinste Tr"agheitsmoment. Kugeln sind jedoch nicht raumfüllend wie
Vorono"izellen. Zur Zeit hat man kein anschauliches Verst"andnis von den
Parametern, die ein niedriges Tr"agheitsmoment bestimmen. In einer neueren Arbeit erzielten Agrell und
Eriksson \cite{agrell98:_optim_lattic_quant} eine Reihe von zum Teil
kontraintuitiven Ergebnissen über Gitterquantizer.
 Sie adaptieren dort einen
Trainingsalgorithmus für Quantizer auf die Gittersituation. Die mit
dieser Methode gefundenen Gitter befinden sich bereits relativ nahe an
der von Conway und Sloane vermuteten theoretischen unteren Schranke für den mittleren quadratischen
Fehler eines Gitterquantizers \cite{conway85}, insbesondere verbessern sie die
bekannten Rekorde  in den Dimensionen $9$ und
$10,$ s.~\cite{bible}, S.~61. In Dimension $10$ fanden sie in $D_{10}^+$ ein  Gitter, dessen sehr
niedriges Tr"agheitsmoment sie mit $0.070814\pm0.000001$ angeben.
 Der Umstand, daß ein relativ
bekanntes Gitter wie  $D_{10}^+$ in dieser
Hinsicht bisher immer übersehen wurde, legt nahe, daß Gitterquantizer
bis dahin nur oberfl"achlich untersucht wurden, ganz 
im Gegensatz zu den klassischen Problemen wie dem Packungs- und dem
Überdeckungsproblem. In Dimension $9$ berechneten sie ein
nichtklassisches Gitter mit der folgenden Erzeugermatrix
\begin{displaymath}
\left(
  \begin{matrix}
    -1&-1&0&0&0&0&0&0&0\\
    1&-1&0&0&0&0&0&0&0\\
    0&1&-1&0&0&0&0&0&0\\
    0&0&1&-1&0&0&0&0&0\\
    0&0&0&1&-1&0&0&0&0\\
    0&0&0&0&1&-1&0&0&0\\
    0&0&0&0&0&1&-1&0&0\\
    0&0&0&0&0&0&1&-1&0\\
    1/2&1/2&1/2&1/2&1/2&1/2&1/2&1/2&0.573\\    
  \end{matrix}
\right)
\end{displaymath}
und einem Tr"agheitsmoment von $0.071626\pm 0.0000002$ (nach
\cite{agrell98:_optim_lattic_quant}), welches sich schon relativ nahe am von
Conway und Sloane vermuteten theoretisch möglichen Minimum $0.070902...$
befindet.  Die ersten acht Zeilen erzeugen das Gitter $D_8.$ Die ersten acht
Koordinaten des neunten Vektors beschreiben ein \glqq tiefes Loch\grqq\footnote{Siehe S. \pageref{deep_hole}.}
  von $D_8.$ Fügt man das deep hole
als Erzeuger zu $D_8$ hinzu, so erh"alt man das Gitter $E_8.$ Die neunte
Koordinate mit dem seltsam anmutenden, numerisch approximierten Wert $0.573$ l"aßt sich
geometrisch nicht einordnen. Das Minimum von $D_8$ ist gleich $2,$ das Gitter
von Agrell und Eriksson besitzt Minimum $(2\cdot0.573)^2 =1.313316.$ Die
Vorono"izelle dieses Gitters ist also relativ flach, w"ahrend man nach der
obigen Motivation eher eine kugelförmigere Vorono"izelle erwartet h"atte.

Mit \texttt{tn} kann  man in verschiedenen Dimensionen   eine
Reihe von Gittern finden mit fast so niedrigen Tr"agheitsmomenten wie
die niedrigsten bekanntesten. Auch in höheren Dimensionen, als sie von
Agrell und Eriksson veröffentlicht wurden, findet man mit \texttt{tn}
Gitter mit Tr"agheitsmomenten in der N"ahe der vermuteten theoretischen unteren
Schranke, wie zum Beispiel $G(A_{11}^{+3})\approx 0.070422\pm 0.000012.$

\bibliographystyle{amsalpha}
\bibliography{lattice}
\backmatter
\end{document}